\documentclass[10pt]{elsarticle}
\usepackage[margin=1.5in]{geometry}
\usepackage{amssymb, amsmath, mathrsfs}
\usepackage{bm}
\usepackage{latexsym}
\usepackage{graphicx}
\usepackage{subfig}
\usepackage{caption}
\usepackage{float}
\usepackage[colorlinks]{hyperref}
\usepackage{tabularx}
\usepackage{multirow}
\usepackage{lineno}
\usepackage{booktabs}
\graphicspath{{./fig/}}
\usepackage{color, xcolor}

\def\D{\mathbb{D}}
\def\I{\mathbb{I}}
\def\S{\mathbb{S}}
\def\d{\mathrm{d}}
\def\x{\bm{x}}

\newtheorem{remark}{Remark}

\newtheorem{theorem}{Theorem}
\newcommand{\jump}[1]{[\![ #1 ]\!]}
\allowdisplaybreaks[4]
\begin{document}

\begin{frontmatter}
    \title{Numerical solutions of an accurate diffuse interface model of the incompressible resistive MHD free surface flow}
    \author{Maojun Li}
    \ead{limj@uestc.edu.cn}
    \author{Jiancheng Wang}
    \ead{202311110403@std.uestc.edu.cn}
    \author{Zeyu Xia\corref{cor1}}
    \ead{zeyuxia@uestc.edu.cn}
    \author{Liwei Xu}
    \ead{xul@uestc.edu.cn}

    \cortext[cor1]{Corresponding author}

    \address{School of Mathematical Sciences, University of Electronic Science and Technology of China, Sichuan, 611731, P.R. China }

    \begin{abstract}
        In this paper, we derive a new model to simulate the incompressible resistive magnetohydrodynamic (MHD) free surface flow. A thermodynamically consistent diffuse interface method is adopted to characterize the moving interface in the modeling process. The formal convergence of the proposed MHD free surface flow model to the sharp interface model is established via a matched asymptotic argument, and the model can be solved without the need for sophisticated free surface capturing schemes. We design a fully decoupled linear finite element scheme that preserves the divergence-free constraint of the magnetic field at a discrete level. The reliability and robustness of the proposed model and algorithm are validated through numerical investigations of the magnetic damping effect on bubble dynamics. In particular, we provide a quantitative numerical comparison of the present results with those obtained from an inductionless MHD model and a sharp interface arbitrary Lagrangian--Eulerian model.
    \end{abstract}


\end{frontmatter}


\section{Introduction}

When a magnetic field is applied to an electrically conducting and non-magnetic fluid (e.g., liquid metals, plasmas, and strong electrolytes), the induced magnetohydrodynamic (MHD) effect can significantly change the fluid dynamics via the Lorentz force, resulting in fundamentally distinct flow behaviors. Given the ubiquitous presence of such a scenario in industrial processes, numerous studies have been carried out to investigate the underlying flow mechanism \cite{2017_Davidson} and develop structure-preserving numerical schemes \cite{2018_Hiptmair, 2017_Hu, 2019_Li, 2025_Mao}. Although MHD flows exhibit extremely complex dynamics in engineering applications, the frequent occurrence of free surfaces further complicates the fluid dynamics. Representative examples include magnetic stirring and damping in the continuous casting and refining of metals \cite{2017_Davidson, 2024_Gou}, MHD instabilities of free surfaces in Hall--H\'{e}roult aluminum reduction cells \cite{2003_Gerbeau, 2019_Herreman}, liquid metal batteries \cite{2019_Herreman, 2018_Tucs}, and plasma-facing components of tokamaks \cite{2019_Lunz, 2004_Morley}.

In this work, we consider incompressible resistive MHD free surface flows in a fixed and bounded domain $\Omega\subset\mathbb{R}^3$, which contains two immiscible fluids occupying two time-dependent open subdomains $\Omega_\pm(t)$, respectively. From a macroscopic perspective, the free surface is a manifold of codimension one \cite{2011_Gross}, and the two subdomains $\Omega_\pm(t)$ are separated by a sharp interface $\Gamma(t)=\overline{\Omega_+(t)}\cap\overline{\Omega_-(t)}$ that should not be in contact with the boundary $\partial\Omega$. Consequently, the flow dynamics can be governed by the following resistive MHD equations in each subdomain $\Omega_\pm(t)$ with the jump conditions on the interface $\Gamma(t)$,
\begin{subequations}
\label{sharp}
    \begin{align}
        \rho_\pm(\partial_t\bm{u}+\bm{u}\cdot\nabla\bm{u}) & =\nabla\cdot\big(2\eta_\pm\D(\bm{u})-p\I\big)+\bm{J}\times\bm{B}+\rho_\pm\bm{g} & & \text{in }\Omega_\pm(t), \label{sharp_1} \\
        \nabla\cdot\bm{u} & =0 & & \text{in }\Omega_\pm(t), \label{sharp_2} \\
        \partial_t\bm{B}+\nabla\times\bm{E} & =\bm{0} & & \text{in }\Omega_\pm(t), \label{sharp_3} \\
        \nabla\times\bm{B} & =\mu_0\bm{J} & &\text{in }\Omega_\pm(t), \label{sharp_4} \\
        \bm{J} & =\sigma_\pm(\bm{E}+\bm{u}\times\bm{B}) & & \text{in }\Omega_\pm(t), \label{sharp_5} \\
        \nabla\cdot\bm{B} & =0,\quad\nabla\cdot\bm{J}=0 & & \text{in }\Omega_\pm(t), \label{sharp_6} \\
        \jump{2\eta\D(\bm{u})-p\I}\bm{n}_\Gamma & =-\lambda\kappa\bm{n}_\Gamma & & \text{on }\Gamma(t), \label{sharp_7} \\
        \jump{\bm{u}} & = {\bm 0} & & \text{on }\Gamma(t), \label{sharp_8} \\
        \jump{\bm{B}} & =\bm{0} & & \text{on }\Gamma(t), \label{sharp_9} \\
        \bm{n}_\Gamma\times\jump{\bm{E}} & =\bm{0} & & \text{on }\Gamma(t), \label{sharp_10} \\
        \jump{\bm{J}}\cdot\bm{n}_\Gamma & =0 & & \text{on }\Gamma(t), \label{sharp_11} \\
        V_\Gamma & =\bm{u}\cdot\bm{n}_\Gamma & & \text{on }\Gamma(t). \label{sharp_12}
    \end{align}
\end{subequations}
In the above system, the unknowns are the fluid velocity $\bm{u}$, the pressure $p$, the magnetic field $\bm{B}$, the electric field $\bm{E}$, and the current density $\bm{J}$. In addition, $\D(\bm{u})=\tfrac{1}{2}\big(\nabla\bm{u}+(\nabla\bm{u})^\top\big)$ and $\I$ stand for the deformation tensor and identity tensor, respectively. $\bm{g}=(0,0,-g)^\top$ represents the gravitational acceleration with $g$ being a positive constant. The material properties $\rho_\pm$, $\eta_\pm$, and $\sigma_\pm$ are the positive constant densities, dynamic viscosities, and electrical conductivities of the two fluids, respectively. Moreover, we denote by $\mu_0$ the permeability of free space, $\lambda$ the surface tension coefficient, $\kappa$ the mean curvature, $\bm{n}_\Gamma$ the unit normal on $\Gamma(t)$ pointing to $\Omega_+(t)$, and $V_\Gamma$ the normal velocity of $\Gamma(t)$. For brevity, we use the notation $\chi=\chi_+\mathscr{X}_{\Omega_+(t)}+\chi_-\mathscr{X}_{\Omega_-(t)}$ throughout this work, where $\chi$ denotes any physical quantity or parameter in system \eqref{sharp}, and $\mathscr{X}_{\Omega_\pm(t)}$ are the characteristic functions corresponding to $\Omega_\pm(t)$. Generally, $\jump{\chi}=\chi_+-\chi_-$ denotes the jump of quantity $\chi$ across $\Gamma(t)$ from $\Omega_+(t)$ to $\Omega_-(t)$. Equations \eqref{sharp_1}-\eqref{sharp_6} are the resistive MHD equations modeling the mutual interaction between the electromagnetic fields and the fluid flows; see \cite{2017_Davidson} for a physical justification. Meanwhile, equations \eqref{sharp_7}-\eqref{sharp_11} are the jump conditions at the interface $\Gamma(t)$, derived from equations \eqref{sharp_1}-\eqref{sharp_6}. Equation \eqref{sharp_12} models the movement of $\Gamma(t)$, corresponding to the immiscibility assumption. Under homogeneous boundary conditions, the sharp interface model \eqref{sharp} satisfies the following stability estimate,
\begin{equation*}
    \frac{\d}{\d{t}}\Bigg(\int_\Omega\frac{1}{2}\rho|\bm{u}|^2\d\x+\int_\Omega\frac{1}{2}\mu_0^{-1}|\bm{B}|^2\d\x+\int_{\Gamma(t)}\lambda\d\bm{s}+\int_\Omega\rho{g}x_3\d\x\Bigg)=-\int_\Omega{2}\eta|\D(\bm{u})|^2\d\x-\int_\Omega\sigma^{-1}|\bm{J}|^2\d\x\le{0},
\end{equation*}
where $x_3$ is the third component of $\x=(x_1,x_2,x_3)$.

Numerical solution of system \eqref{sharp} is challenging primarily due to the presence of the moving interface $\Gamma(t)$ and the associated jump conditions. To characterize the evolving interfaces in MHD free surface problems, substantial numerical methodologies have been developed in the past decades. Morley et al. \cite{2004_Morley} employed the level set technique to model the free surface of fusion liquid walls in tokamaks, where the interface was represented by a signed distance function, and the material properties were obtained through an approximate Heaviside function. In addition, the re-initialization method and the continuum surface force (CSF) model were exploited to conduct stable numerical simulations, and interested readers are referred  to \cite{2011_Gross} for more details on these methods. Similar techniques were also utilized by Munger and Vincent \cite{2006_Munger}, where a finite volume level set approach was developed to simulate the MHD instability in aluminum reduction cells while neglecting the surface tension force. Flueck et al. \cite{2010_Flueck} utilized the level set method as an indicator and constructed an interface-fitted mesh using the zero level set of the solution, to simulate the instability of the interface in aluminum reduction cells. Cappanera et al. \cite{2018_Cappanera} adopted a conservative level set method, which utilized an approximate Heaviside function rather than the signed distance function, to investigate the metal pad roll instability in stably stratified liquid layers; see also \cite{2019_Herreman}. Furthermore, Zhang and Ni \cite{2014_Zhang} applied the volume of fluid method, representing the interface by volume fractions, to explore the bubble dynamics in metallurgical processes. In their approach, the material properties and surface tension force were treated smoothly, and an adaptive mesh refinement technique combined with a consistent and conservative scheme for the Lorentz force was developed to accurately resolve the fluid dynamics. They also designed an algorithm to address the phase transition problem using cut cell and ghost cell techniques in \cite{2018_Zhang}. In addition to implicit interface capturing approaches such as the level set and volume of fluid methods, explicit interface tracking and moving mesh techniques have also been extensively developed in the simulation of MHD problems. Gerbeau et al. \cite{2003_Gerbeau} employed an arbitrary Lagrangian--Eulerian (ALE) approach for a two-fluid MHD problem to simulate MHD instability of aluminum reduction cells; see also Ba\v{n}as and Prohl \cite{2010_Lubomir} for a convergent bound-preserving scheme for this system. Samulyak et al. \cite{2007_Samulyak} utilized a front tracking technique to capture the interface in a compressible inviscid MHD model. Pan et al. \cite{2018_Pan} exploited the immersed boundary method to accurately solve the inductionless MHD equations with a moving boundary. While these approaches are widely used, they typically require specialized numerical techniques to ensure high interface resolution or numerical stability.

In addition to the aforementioned approaches, the diffuse interface method \cite{1998_Anderson} emerges as another effective tool to characterize the free surface between binary immiscible fluids. In this methodology, the sharp interface $\Gamma(t)$ is represented by a significantly thin transition layer where the fluid properties vary smoothly, and the singular surface tension effects can be modeled by either the CSF model or the smooth Korteweg stress tensor. Hence, a stable simulation can be readily conducted without the additional numerical effort that is usually required by the aforementioned methods \cite{2011_Gross}. Regarding the diffuse interface method for solving the MHD free surface flow, \cite{2019_Yang} proposed a scheme to simulate the resistive MHD problem with matched densities, and established the existence of weak solutions. Subsequently, many studies have been devoted to the convergence analysis of numerical schemes for this model; see, e.g., \cite{2024_WangC, 2025_Yang}. In addition, \cite{2020_Chen} proposed a diffuse interface method for an inductionless MHD problem with matched densities. We point out that the two-phase MHD model employed in the above works is a direct coupling of MHD equations with a phase model, and its relationship with the classical sharp interface formulation has not been quantified.

In this paper, we employ a thermodynamically consistent diffuse interface method \cite{2012_Abels} to numerically solve the MHD free surface flow. Although the sharp interface $\Gamma(t)$ is replaced by a thin transition layer, we demonstrate that system \eqref{sharp} is formally recovered as the thickness of the transition layer tends to zero, by utilizing the method of matched asymptotic expansions \cite{2012_Abels, 2018_Xu}. Then, we develop a fully decoupled linear finite element scheme that involves only Poisson-type equations and time-independent saddle point systems. Accordingly, algebraic multigrid and additive Schwarz methods can be exploited as preconditioners for the resulting linear systems, in turn facilitating large-scale computations. Furthermore, we apply a potential method \cite{2020_Li} to solve the pre-Maxwell equations \eqref{sharp_3}-\eqref{sharp_6}, ensuring an exactly solenoidal discrete magnetic field. Finally, we present numerical examples of two- and three-dimensional single rising bubble benchmarks, together with a quantitative comparison of our results with those obtained from an inductionless MHD model and a sharp interface ALE model. The results demonstrate the capability and robustness of the proposed method in capturing complex MHD effects in free surface flows.

The remainder of this paper is organized as follows. In Section \ref{section_2}, we present the numerical framework and demonstrate its convergence to the sharp interface model. The discrete scheme is also detailed. In Section \ref{section_3}, several numerical results on the magnetic damping effect are presented, including the effect of the interface thickness and the comparison among different models. Conclusions and future research directions are discussed in Section \ref{section_4}.

\section{Numerical framework}
\label{section_2}

\subsection{Governing equations}

With the sharp interface $\Gamma(t)$ replaced by a thin transition layer $\Gamma_\epsilon(t)$ where the two fluids are mixed, model \eqref{sharp} is approximated by the following diffuse interface system (the derivation of these equations is outlined in \ref{derivation}):
\begin{subequations}
\label{phase}
    \begin{align}
        \partial_t\phi+\nabla\cdot(\phi\bm{u}) & =\nabla\cdot\big(m_\epsilon(\phi)\nabla\mu\big), \label{phase_1} \\
        \mu & =-\widehat{\lambda}\epsilon\Delta\phi+\widehat{\lambda}\epsilon^{-1}(\phi^3-\phi), \label{phase_2} \\
        \rho(\phi)\partial_t\bm{u}+\big(\rho(\phi)\bm{u}-\rho_dm_\epsilon(\phi)\nabla\mu\big)\cdot\nabla\bm{u} & =\nabla\cdot\big(2\eta(\phi)\D(\bm{u})-p\I\big)+\mu\nabla\phi+\bm{J}\times\bm{B}+\rho(\phi)\bm{g}, \label{phase_3} \\
        \nabla\cdot\bm{u} & =0, \label{phase_4} \\
        \partial_t\bm{B}+\nabla\times\bm{E} & =\bm{0}, \label{phase_5} \\
        \nabla\times\bm{B} & =\mu_0\bm{J}, \label{phase_6} \\
        \bm{J} & =\sigma(\phi)(\bm{E}+\bm{u}\times\bm{B}), \label{phase_7} \\
        \nabla\cdot\bm{B} & =0,\quad\nabla\cdot\bm{J}=0. \label{phase_8}
    \end{align}
\end{subequations}
In the above equations, $\phi$ is the order parameter that labels the two fluids such that
\begin{equation*}
    \phi(\x,t)=
    \begin{cases}
        1,  & \x\in\Omega_+(t),\\
        -1, & \x\in\Omega_-(t),
    \end{cases}
    \qquad    \mbox{and} \qquad
    |\phi(\x,t)|<1\;\;\text{in}\;\;\Gamma_\epsilon(t),
\end{equation*}
$\mu$ is the chemical potential serving as an auxiliary variable, $\epsilon$ is the thickness of $\Gamma_\epsilon(t)$, $m_\epsilon(\phi)>0$ represents the mobility which plays the role of diffusivity for $\Gamma_\epsilon(t)$, and $\widehat{\lambda}=\tfrac{3}{2\sqrt{2}}\lambda$ stands for the scaled surface tension coefficient. Moreover, the density $\rho$, the dynamic viscosity $\eta$, and the electrical conductivity $\sigma$ in the mixture are all linearly dependent on $\phi$ as illustrated below for $\rho$,
\begin{equation*}
    \rho(\phi)=\frac{\rho_+-\rho_-}{2}\phi+\frac{\rho_++\rho_-}{2}.
\end{equation*}
We define $\rho_d=\tfrac{\partial\rho}{\partial\phi}=\tfrac{\rho_+-\rho_-}{2}$ and the extra mass flux term $\rho_dm_\epsilon\nabla\mu$ accounts for thermodynamic consistency. $\mu\nabla\phi$ is a simplified Korteweg stress modeling the surface tension force. Equations \eqref{phase_1}-\eqref{phase_2} represent the widely recognized Cahn--Hilliard equation which describes the process of phase separation.

To close the system, we impose some appropriate initial and boundary conditions. The initial data are specified as
\begin{equation*}
    \bm{u}(0)=\bm{u}^0,\qquad\bm{B}(0)=\bm{B}^0,
    \qquad    \mbox{and} \qquad
    \phi(0)=\phi^0,
\end{equation*}
with $\nabla\cdot\bm{u}^0=0$ and $\nabla\cdot\bm{B}^0=0$. For simplicity, we apply the following homogeneous boundary conditions
\begin{equation*}
    \begin{aligned}
        \bm{u}=\bm{0}, \qquad\bm{n}\times\bm{E}=\bm{0}, \qquad\nabla\phi\cdot\bm{n}=0,
        \qquad    \mbox{and} \qquad
        \nabla\mu\cdot\bm{n}=0
    \end{aligned}
\end{equation*}
on $\partial\Omega$, where $\bm{n}$ is the unit outward normal to $\Omega$. Boundary conditions adapted to the physical setting of benchmarks will be specified in the numerical experiments.

\begin{remark}
    System \eqref{phase} satisfies the following energy balance law, whose detailed derivation is provided in \ref{law}:
    \begin{equation*}
        \begin{aligned}
            &\frac{\d}{\d{t}}\Bigg(\int_\Omega\frac{1}{2}\rho|\bm{u}|^2\d\x+\int_\Omega\frac{1}{2}\mu_0^{-1}|\bm{B}|^2\d\x+\int_\Omega{f}(\phi,\nabla\phi)\d\x\Bigg) \\
            &=-\int_\Omega{m}_\epsilon|\nabla\mu|^2\d\x-\int_\Omega{2}\eta|\D(\bm{u})|^2\d\x-\int_\Omega\sigma^{-1}|\bm{J}|^2\d\x+\int_\Omega\rho\bm{g}\cdot\bm{u}\d\x,
        \end{aligned}
    \end{equation*}
    where $f(\phi,\nabla\phi)=\tfrac{\widehat{\lambda}\epsilon}{2}|\nabla\phi|^2+\tfrac{\widehat{\lambda}\epsilon^{-1}}{4}(\phi^2-1)^2$. Because of the diffusive mass flux inherent in the volume-averaged formulation (a numerical demonstration of the influence of this flux can be found in \cite{2014_Grun}), the gravitational contribution cannot in general be absorbed entirely into a standard gravitational potential energy term. In simulations, gravity always plays an important role, and the energy of system \eqref{phase} is no longer dissipative, unlike that of model \eqref{sharp}. In contrast, interested readers are referred to \cite{1998_Lowengrub} for another thermodynamically consistent diffuse interface method leading to a much more complex quasi-incompressible system that preserves a conservative gravitational force.
\end{remark}

Since the thickness of the transition layer $\Gamma_\epsilon(t)$ is determined by parameter $\epsilon$, physically correct behavior is expected in simulations for sufficiently small $\epsilon$. In this paper, we take a constant mobility with an established scaling law $m_\epsilon=\mathcal{O}(\epsilon^2)$ \cite{2013_Magaletti, 2024_Wang, 2018_Xu}, and hence we are able to establish the formal convergence of system \eqref{phase} to model \eqref{sharp} as $\epsilon\to{0}$.

\begin{theorem}
    \label{sharp_asymptotics}
    Assuming $m_\epsilon=\mathcal{O}(\epsilon^2)$, the formal asymptotic limit of system \eqref{phase} for $\epsilon\to{0}$ is model \eqref{sharp}.
\end{theorem}

\subsection{Formal derivation of the sharp interface limit}

We exploit the method of formally matched asymptotic expansions to identify the sharp interface limit. This argument depends on two critical assumptions: (i) under vanishing interface thickness the domain can be decomposed into two distinct subdomains separated by an interface, each containing only one fluid; and (ii) the solutions admit separate asymptotic expansions in terms of the interface thickness in the bulk and near the interface, which must be matched in the overlapping region.

\subsubsection{Preliminaries}

Letting $(\phi_\epsilon,\mu_\epsilon,\bm{u}_\epsilon,p_\epsilon,\bm{B}_\epsilon,\bm{E}_\epsilon,\bm{J}_\epsilon)$ be a solution to equations \eqref{phase}, we assume that it converges formally to a limit $(\phi,\mu,\bm{u},p,\bm{B},\bm{E},\bm{J})$ as $\epsilon\to{0}$, and further,
\begin{equation*}
    \Omega_+(t)=\{\x\in\Omega:\phi_\epsilon(\x)>0\}
    \qquad \mbox{and} \qquad
    \Omega_-(t)=\{\x\in\Omega:\phi_\epsilon(\x)<0\}.
\end{equation*}
Then, the corresponding interface separating $\Omega_\pm(t)$ is denoted by
\begin{equation*}
    \Gamma(t)=\{\x\in\Omega:\phi_\epsilon(\x)=0\}.
\end{equation*}
Here, we reuse for simplicity the notations $\Omega_\pm(t)$ and $\Gamma(t)$.

In outer regions far from $\Gamma(t)$, we assume the following outer expansions with respect to the interface thickness $\epsilon$ of $\phi_\epsilon$, $\bm{u}_\epsilon$, $p_\epsilon$, $\bm{B}_\epsilon$, $\bm{E}_\epsilon$, $\bm{J}_\epsilon$; for example, for $\phi_\epsilon$,
\begin{equation*}
    \phi_\epsilon=\phi_0^\pm+\epsilon\phi_1^\pm+\epsilon^2\phi_2^\pm+\cdots.
\end{equation*}
As an auxiliary variable, we have evidently for $\mu_\epsilon$,
\begin{equation*}
    \mu_\epsilon=\epsilon^{-1}\mu_0^\pm+\mu_1^\pm+\epsilon\mu_2^\pm+\cdots
\end{equation*}
with $\mu_0^\pm=\widehat{\lambda}\big((\phi_0^\pm)^3-\phi_0^\pm\big)$. In the transition layer near $\Gamma(t)$, we introduce a new stretched coordinate to scale the variables. Let $d(\x)$ be the signed distance function to $\Gamma(t)$ such that $d(\x)>0$ if $\x\in\Omega_+(t)$. Then we have $\bm{n}_\Gamma=\nabla{d}$ and $V_\Gamma=-\partial_td$. Defining the scaled distance $z=\tfrac{d}{\epsilon}$, the following change of variables holds \cite{2012_Abels}:
\begin{subequations}
\label{change_variable}
    \begin{align}
        \partial_t\chi & =-\epsilon^{-1}V_\Gamma\partial_z\chi+\mathcal{O}(1), \\
        \nabla\chi & =\epsilon^{-1}\partial_z\chi\bm{n}_\Gamma+\nabla_\Gamma\chi+\mathcal{O}(\epsilon), \\
        \nabla\cdot\bm{\chi} & =\epsilon^{-1}\partial_z\bm{\chi}\cdot\bm{n}_\Gamma+\nabla_\Gamma\cdot\bm{\chi}+\mathcal{O}(\epsilon), \\
        \Delta\chi & =\epsilon^{-2}\partial_{zz}\chi-\epsilon^{-1}\kappa\partial_z\chi+\mathcal{O}(1),
    \end{align}
\end{subequations}
where $\chi$ and $\bm{\chi}$ represent the generic scalar- and vector-valued variables, respectively. Then in this region we have the inner expansion in terms of the interface thickness $\epsilon$ for $\phi_\epsilon$, $\bm{u}_\epsilon$, $p_\epsilon$, $\bm{B}_\epsilon$, $\bm{E}_\epsilon$, $\bm{J}_\epsilon$ in the form (also exemplarily for $\phi_\epsilon$)
\begin{equation*}
    \phi_\epsilon=\phi_0^i+\epsilon\phi_1^i+\epsilon^2\phi_2^i+\cdots,
\end{equation*}
with $\mu_\epsilon=\epsilon^{-1}\mu_0^i+\mu_1^i+\epsilon\mu_2^i+\cdots$, where
\begin{equation*}
    \mu_0^i=\widehat{\lambda}\big(-\partial_{zz}\phi_0^i+(\phi_0^i)^3-\phi_0^i\big)
    \qquad \mbox{and} \qquad
    \mu_1^i=\widehat{\lambda}\big(-\partial_{zz}\phi_1^i+\kappa\partial_z\phi_0^i+(3(\phi_0^i)^2-1)\phi_1^i\big).
\end{equation*}
In the overlapping region, the outer and inner expansions must satisfy the matching conditions \cite{2012_Abels}
\begin{subequations}
    \begin{align}
        \lim_{z\to\pm\infty}\chi_0^i & =\chi_0^\pm,\label{match_1} \\
        \lim_{z\to\pm\infty}\partial_z\chi_0^i & =0,\label{match_2} \\
        \lim_{z\to\pm\infty}\partial_z\chi_1^i & =\nabla\chi_0^\pm\cdot\bm{n}_\Gamma.\label{match_3}
    \end{align}
\end{subequations}

Next, we will establish that the limit solution $(\phi,\mu,\bm{u},p,\bm{B},\bm{E},\bm{J})$ exactly satisfies equations \eqref{sharp} by following the procedure in \cite{2012_Abels, 2018_Xu}.

\subsubsection{Bulk equations}

\textbf{Expansion of \eqref{phase_3} at order $\mathcal{O}(\epsilon^{-1})$:} Recalling the assumption $m_\epsilon:=m\epsilon^2$, we naturally obtain that the leading order of the flux $m_\epsilon\nabla\mu$ is $\mathcal{O}(\epsilon)$. Hence, we arrive at
\begin{equation*}
    \mu_0^\pm\nabla\phi_0^\pm=0,
\end{equation*}
which is equivalent to
\begin{equation*}
    \nabla\frac{\big((\phi_0^\pm)^2-1\big)^2}{4}=\bm{0}.
\end{equation*}
This implies that $\phi_0^\pm=c_\pm$ where $c_\pm$ are two constants such that $c_+>0$ and $c_-<0$.

\textbf{Expansion of \eqref{phase_3}-\eqref{phase_8} at order $\mathcal{O}(1)$:} Using the fact that $m_\epsilon\nabla\mu$ is of order $\mathcal{O}(\epsilon)$, we conclude
\begin{equation*}
    \begin{aligned}
        \rho_0^\pm\big(\partial_t\bm{u}_0^\pm+\bm{u}_0^\pm\cdot\nabla\bm{u}_0^\pm\big) & =\nabla\cdot\big(2\eta_0^\pm\D(\bm{u}_0^\pm)-p_0^\pm\I\big)+\mu_0^\pm\nabla\phi_1^\pm+\bm{J}_0^\pm\times\bm{B}_0^\pm+\rho_0^\pm\bm{g}, \\
        \nabla\cdot\bm{u}_0^\pm & =0, \\
        \partial_t\bm{B}_0^\pm+\nabla\times\bm{E}_0^\pm & =0, \\
        \nabla\times\bm{B}_0^\pm & =\mu_0\bm{J}_0^\pm, \\
        \bm{J}_0^\pm & =\sigma_0^\pm\big(\bm{E}_0^\pm+\bm{u}_0^\pm\times\bm{B}_0^\pm\big), \\
        \nabla\cdot\bm{B}_0^\pm & =0,\quad\nabla\cdot\bm{J}_0^\pm=0,
    \end{aligned}
\end{equation*}
where we abbreviate $\chi_0^\pm=\tfrac{(\chi_+-\chi_-)\phi_0^\pm}2+\tfrac{\chi_++\chi_-}2$ for $\chi\in\{\rho, \eta, \sigma\}$.

\subsubsection{Jump conditions}

\textbf{Expansion of \eqref{phase_4} at order $\mathcal{O}(\epsilon^{-1})$:} Invoking change of variables \eqref{change_variable}, we get
\begin{equation}
\label{inner_1}
    \partial_z\bm{u}_0^i\cdot\bm{n}_\Gamma=0,
\end{equation}
which implies that $\bm{u}_0^i\cdot\bm{n}_\Gamma$ is independent of $z$.

\textbf{Expansion of \eqref{phase_1} at order $\mathcal{O}(\epsilon^{-1})$:} With a similar procedure, we obtain
\begin{equation*}
    -V_\Gamma\partial_z\phi_0^i+\partial_z\big(\phi_0^i\bm{u}_0^i\big)\cdot\bm{n}_\Gamma=m\partial_{zz}\mu_0^i.
\end{equation*}
Integrating this equation with respect to $z$ from $-\infty$ to $\infty$ yields
\begin{equation*}
    (-V_\Gamma+\bm{u}_0^i\cdot\bm{n}_\Gamma)\int_{-\infty}^\infty\partial_z\phi_0^i\d{z}=m\int_{-\infty}^\infty\partial_{zz}\mu_0^i\d{z}.
\end{equation*}
Then matching conditions \eqref{match_1}-\eqref{match_2} require
\begin{equation*}
    \big(-V_\Gamma+\bm{u}_0^i\cdot\bm{n}_\Gamma\big)\big(\phi_0^+-\phi_0^-\big)=0,
\end{equation*}
in turn giving
\begin{equation*}
    V_\Gamma=\bm{u}_0^i\cdot\bm{n}_\Gamma.
\end{equation*}
Additionally, we have
\begin{equation*}
    \partial_{zz}\mu_0^i=0,
\end{equation*}
or equivalently, $\mu_0^i$ linearly dependent on $z$. Matching condition \eqref{match_1} then requires $\mu_0^i$ to be a constant.

\textbf{Expansion of \eqref{phase_3} at order $\mathcal{O}(\epsilon^{-2})$:} Defining $\mathcal{D}(\bm{A})=\tfrac{1}{2}\big(\bm{A}+\bm{A}^\top\big)$, we obtain after the change of variables \eqref{change_variable},
\begin{equation*}
    \nabla\cdot\big(2\eta\D(\bm{u})\big)=\epsilon^{-2}\partial_z\big(2\eta\mathcal{D}(\partial_z\bm{u}\otimes\bm{n}_\Gamma)\bm{n}_\Gamma\big)+\epsilon^{-1}\partial_z\big(2\eta\mathcal{D}(\nabla_\Gamma\bm{u})\bm{n}_\Gamma\big)+\epsilon^{-1}\nabla_\Gamma\cdot\big(2\eta\mathcal{D}(\partial_z\bm{u}\otimes\bm{n}_\Gamma\big)+\mathcal{O}(1).
\end{equation*}
Noting that $m_\epsilon\nabla\mu$ is of order $\mathcal{O}(\epsilon)$ in the inner region since $\mu_0^i$ is independent of $z$, we conclude
\begin{equation*}
    \partial_z\big(\eta_0^i\partial_z\bm{u}_0^i\big)+\mu_0^i\partial_z\phi_0^i\bm{n}_\Gamma=0,
\end{equation*}
where $\eta_0^i=\tfrac{(\eta_+-\eta_-)\phi_0^i}2+\tfrac{\eta_++\eta_-}2$ and we have used
\begin{equation*}
    \big(\bm{n}_\Gamma\otimes\partial_z\bm{u}_0^i\big)\bm{n}_\Gamma=\big(\partial_z\bm{u}_0^i\cdot\bm{n}_\Gamma\big)\bm{n}_\Gamma=\bm{0}.
\end{equation*}
By integrating this equation with respect to $z$, we can obtain
\begin{equation*}
    \eta_0^i\partial_z\bm{u}_0^i\big|_{-\infty}^\infty+\mu_0^i\phi_0^i\big|_{-\infty}^\infty\bm{n}_\Gamma=\bm{0}.
\end{equation*}
Taking the inner product with $\bm{n}_\Gamma$ and using equation \eqref{inner_1} together with matching condition \eqref{match_1}, we obtain
\begin{equation*}
    \mu_0^i\big(\phi_0^+-\phi_0^-\big)=0,
\end{equation*}
in turn indicating $\mu_0^i=0$, or equivalently
\begin{equation*}
    -\partial_{zz}\phi_0^i+\big(\phi_0^i\big)^3-\phi_0^i=0.
\end{equation*}
The solvability of this equation \cite{2018_Xu} leads to
\begin{equation*}
    \phi_0^\pm=\lim_{z\to\pm\infty}\phi_0^i=\pm{1}\;\;\text{and }\;\;\phi_0^i=\tanh\frac{z}{\sqrt{2}}.
\end{equation*}
Then we have a second-order ordinary differential equation
\begin{equation*}
    \partial_z\big(\eta_0^i\partial_z\bm{u}_0^i\big)=\bm{0},
\end{equation*}
which, by matching conditions \eqref{match_1}-\eqref{match_2} and the positivity of $\eta$, implies that $\bm{u}_0^i$ is independent of $z$. Hence, we obtain
\begin{equation*}
    \jump{\bm{u}_0^\pm}=\bm{0}.
\end{equation*}

\textbf{Expansion of \eqref{phase_3} at order $\mathcal{O}(\epsilon^{-1})$:} The combination of matching conditions \eqref{match_1} and \eqref{match_3} yields
\begin{equation}
\label{match_4}
    \lim_{z\to\pm\infty}\big(\partial_z\bm{u}_1^i\otimes\bm{n}_\Gamma+\nabla_\Gamma\bm{u}_0^i\big)=\nabla\bm{u}_0^\pm.
\end{equation}
Noting that $\bm{u}_0^i$ and $\mu_0^i$ are independent of $z$, we can derive
\begin{equation*}
    \partial_z\big(2\eta_0^i\mathcal{D}(\partial_z\bm{u}_1^i\otimes\bm{n}_\Gamma)\bm{n}_\Gamma\big)+\partial_z\big(2\eta_0^i\mathcal{D}(\nabla_\Gamma\bm{u}_0^i)\bm{n}_\Gamma\big)-\partial_z{p}_0^i\bm{n}_\Gamma+\mu_1^i\partial_z\phi_0^i\bm{n}_\Gamma=\bm{0}.
\end{equation*}
Then integrating with respect to $z$ and using matching conditions \eqref{match_1} and \eqref{match_4} yields
\begin{equation*}
    \jump{2\eta\D(\bm{u}_0^\pm)-p_0^\pm\I}\bm{n}_\Gamma=-\int_{-\infty}^\infty\mu_1^i\partial_z\phi_0^i\bm{n}_\Gamma\d{z}.
\end{equation*}
For the right-hand side of this equation, integration by parts with matching conditions \eqref{match_1}-\eqref{match_2} leads to
\begin{equation*}
    \begin{aligned}
        \int_{-\infty}^\infty\mu_1^i\partial_z\phi_0^i\d{z} & =\widehat{\lambda}\int_{-\infty}^\infty\big(-\partial_{zz}\phi_1^i+\kappa\partial_z\phi_0^i+(3(\phi_0^i)^2-1)\phi_1^i\big)\partial_z\phi_0^i\d{z}, \\
        & =\widehat{\lambda}\kappa\int_{-\infty}^\infty(\partial_z\phi_0^i)^2\d{z}+\widehat{\lambda}\int_{-\infty}^\infty\big(-\partial_{zz}\phi_1^i+(3(\phi_0^i)^2-1)\phi_1^i\big)\partial_z\phi_0^i\d{z} \\
        & =\frac{2\sqrt{2}\widehat{\lambda}\kappa}{3}+\widehat{\lambda}\int_{-\infty}^\infty\big(\partial_{zz}\phi_0^i-(\phi_0^i)^3+\phi_0^i\big)\partial_z\phi_1^i\d{z} \\
        & =\lambda\kappa,
    \end{aligned}
\end{equation*}
where
\begin{equation*}
    \int_{-\infty}^\infty(\partial_z\phi_0^i)^2\d{z}=\frac{1}{2}\int_{-\infty}^\infty\mathrm{sech}^4\frac{z}{\sqrt{2}}\d{z}=\frac{2\sqrt{2}}{3}.
\end{equation*}
Then we conclude
\begin{equation*}
    \jump{2\eta\D(\bm{u}_0^\pm)-p_0^\pm\I}\bm{n}_\Gamma=-\lambda\kappa\bm{n}_\Gamma.
\end{equation*}

\textbf{Expansion of \eqref{phase_5}-\eqref{phase_6} and \eqref{phase_8} at order $\mathcal{O}(\epsilon^{-1})$:} A direct change of variables \eqref{change_variable} gives
\begin{equation*}
    \begin{aligned}
        -V_\Gamma\partial_z\bm{B}_0^i+\bm{n}_\Gamma\times\partial_z\bm{E}_0^i & =\bm{0}, \\
        \bm{n}_\Gamma\times\partial_z\bm{B}_0^i & =\bm{0}, \\
        \partial_z\bm{B}_0^i\cdot\bm{n}_\Gamma & =0, \\
        \partial_z\bm{J}_0^i\cdot\bm{n}_\Gamma & =0.
    \end{aligned}
\end{equation*}
Integrating with respect to $z$ with matching condition \eqref{match_1}, we obtain
\begin{subequations}
    \begin{align}
        \jump{\bm{B}_0^\pm} & =\bm{0}, \\
        \bm{n}_\Gamma\times\jump{\bm{E}_0^\pm} & =\bm{0}, \\
        \jump{\bm{J}_0^\pm}\cdot\bm{n}_\Gamma & =0.
    \end{align}
\end{subequations}

\subsubsection{Summary of the sharp interface limit}

We summarize the main results as follows. For the region far from $\Gamma(t)$, we have $\phi_0^\pm=\pm{1}$ and $\mu_0^\pm=0$. Then, we obtain the following resistive MHD equations:
\begin{equation*}
    \begin{aligned}
        \rho_\pm\big(\partial_t\bm{u}_0^\pm+\bm{u}_0^\pm\cdot\nabla\bm{u}_0^\pm\big) & =\nabla\cdot\big(2\eta_0^\pm\D(\bm{u}_0^\pm)-p_0^\pm\I\big)+\bm{J}_0^\pm\times\bm{B}_0^\pm+\rho_\pm\bm{g}, \\
        \nabla\cdot\bm{u}_0^\pm & =0, \\
        \partial_t\bm{B}_0^\pm+\nabla\times\bm{E}_0^\pm & =\bm{0}, \\
        \nabla\times\bm{B}_0^\pm & =\mu_0\bm{J}_0^\pm, \\
        \bm{J}_0^\pm & =\sigma_\pm\big(\bm{E}_0^\pm+\bm{u}_0^\pm\times\bm{B}_0^\pm\big), \\
        \nabla\cdot\bm{B}_0^\pm & = 0,\quad\nabla\cdot\bm{J}_0^\pm=0,
    \end{aligned}
\end{equation*}
and the jump conditions:
\begin{subequations}
    \begin{align*}
        \jump{2\eta\D(\bm{u}_0^\pm)-p_0^\pm\I}\bm{n}_\Gamma & =-\lambda\kappa\bm{n}_\Gamma, \\
        \jump{\bm{u}_0^\pm} & =\bm{0}, \\
        \jump{\bm{B}_0^\pm} & =\bm{0}, \\
        \bm{n}_\Gamma\times\jump{\bm{E}_0^\pm} & =\bm{0}, \\
        \jump{\bm{J}_0^\pm}\cdot\bm{n}_\Gamma & =0, \\
        V_\Gamma & =\bm{u}_0^i\cdot\bm{n}_\Gamma.
    \end{align*}
\end{subequations}
We can deduce that the limit system is exactly model \eqref{sharp}.

\subsection{Discrete scheme}

The incompressible resistive MHD equations involve several fundamental constraints, including the incompressibility of the fluid velocity $\bm{u}$ and the divergence-free conditions on the magnetic field $\bm{B}$ and the current density $\bm{J}$. Some structure-preserving techniques for the incompressible resistive MHD equations, addressing these constraints at the discrete level, can be found in \cite{2018_Hiptmair, 2017_Hu, 2025_Mao}. In the present work, with the objective of developing a fully decoupled linear scheme for the proposed model, we restrict our attention to the discrete preservation of the divergence-free constraint on the magnetic field $\bm{B}$, which is particularly important for MHD simulations. To this end, we introduce the vector potential $\bm{A}$ under the Weyl gauge, so that $\bm{B}=\nabla\times\bm{A}$ and $\bm{E}=-\partial_t\bm{A}$, and obtain \cite{2020_Li}
\begin{equation*}
    \begin{aligned}
        \partial_t\bm{A}-\bm{u}\times\nabla\times\bm{A}+\mu_0^{-1}\sigma^{-1}(\phi)\nabla\times\nabla\times\bm{A} & =0, \\
        \bm{B} & =\nabla\times\bm{A},
    \end{aligned}
\end{equation*}
The boundary condition $\bm{n}\times\bm{A}=\bm{0}$ is imposed consistently with $\bm{n}\times\bm{E}=\bm{0}$ under the Weyl gauge. Meanwhile, consistent with the initial magnetic field $\bm{B}^0$, the initial vector potential is specified as $\bm{A}(0)=\bm{A}^0$, where $\bm{A}^0$ satisfies $\nabla\times\bm{A}^0=\bm{B}^0$.
Using Amp\`{e}re's law, we then rewrite the Lorentz force in terms of the Maxwell stress tensor \cite{2017_Davidson}:
\begin{equation*}
    \bm{J}\times\bm{B}=\mu_0^{-1}\nabla\cdot\Bigg(\bm{B}\otimes\bm{B}-\frac{1}{2}|\bm{B}|^2\I\Bigg).
\end{equation*}
Therefore, the current density $\bm{J}$ does not need to be treated as an additional unknown in the resulting formulation. For simplicity, we absorb the magnetic pressure contribution into the pressure variable by redefining $p\leftarrow{p}+\tfrac{1}{2}\mu_0^{-1}|\bm{B}|^2$. Consequently, only the term $\mu_0^{-1}\nabla\cdot(\bm{B}\otimes\bm{B})$ arising from the Maxwell stress tensor appears explicitly in the momentum equation.

Now let $L_r$, $u_r$, $B_r$, $\rho_r$, $\eta_r$, $\sigma_r$ denote the characteristic quantities of length, velocity, magnetic field, density, dynamic viscosity, and electrical conductivity, respectively. We introduce the following scalings
\begin{equation*}
    \begin{aligned}
         & \x\leftarrow\frac{\x}{L_r},
         &&{t}\leftarrow\frac{tu_r}{L_r},
         &&&\bm{u}\leftarrow\frac{\bm{u}}{u_r},
         &&&&\rho\leftarrow\frac{\rho}{\rho_r},
         &&&&&\eta\leftarrow\frac{\eta}{\eta_r}, \\
         & {p}\leftarrow\frac{p}{\rho_ru_r^2},
         &&\mu\leftarrow\frac{\mu\epsilon}{\widehat{\lambda}},
         &&&\bm{B}\leftarrow\frac{\bm{B}}{B_r},
         &&&&\bm{A}\leftarrow\frac{\bm{A}}{L_rB_r},
         &&&&&\sigma\leftarrow\frac{\sigma}{\sigma_r}.
    \end{aligned}
\end{equation*}
Certainly, the dimensionless system is given by
\begin{subequations}
    \begin{align*}
        \partial_t\phi+\nabla\cdot(\phi\bm{u}) & =\frac{1}{\mathrm{Pe}}\Delta\mu, \\
        \mu & =-\mathrm{Cn}^2\Delta\phi+\phi^3-\phi, \\
        \rho\partial_t\bm{u}+\Bigg(\rho\bm{u}-\rho_d\frac{1}{\mathrm{Pe}}\nabla\mu\Bigg)\cdot\nabla\bm{u} & =\nabla\cdot\Bigg(\frac{2}{\mathrm{Re}}\eta\D(\bm{u})-p\I\Bigg)+\frac{1}{\mathrm{We}\mathrm{Cn}}\mu\nabla\phi+\frac{\mathrm{N}}{\mathrm{Rm}}\nabla\cdot(\bm{B}\otimes\bm{B})-\frac{1}{\mathrm{Fr}}\rho\bm{e}_3, \\
        \nabla\cdot\bm{u} & =0, \\
        \partial_t\bm{A}-\bm{u}\times\nabla\times\bm{A} & =-\frac{1}{\mathrm{Rm}}\sigma^{-1}\nabla\times\nabla\times\bm{A}, \\
        \bm{B} & =\nabla\times\bm{A},
    \end{align*}
\end{subequations}
where $\bm{e}_3=(0,0,1)^\top$. Here, the dimensionless numbers include the Cahn number $\mathrm{Cn}=\tfrac{\epsilon}L_r$, P\'{e}clet number $\mathrm{Pe}=\tfrac{\epsilon{L}_ru_r}{\widehat{\lambda}m_\epsilon}$, Reynolds number $\mathrm{Re}=\tfrac{L_r\rho_ru_r}{\eta_r}$, Weber number $\mathrm{We}=\tfrac{L_r\rho_ru_r^2}{\widehat{\lambda}}$, Stuart number $\mathrm{N}=\tfrac{L_r\sigma_rB_r^2}{\rho_ru_r}$, Froude number $\mathrm{Fr}=\tfrac{u_r^2}{gL_r}$, and magnetic Reynolds number $\mathrm{Rm}=L_ru_r\mu_0\sigma_r$.

Let $L^2(\Omega)$ be the Lebesgue space of square integrable functions equipped with inner product $(\cdot,\cdot)$, and let $H^1(\Omega)$, $\bm{H}(\mathrm{div};\Omega)$, and $\bm{H}(\mathrm{curl};\Omega)$ be subspaces of $L^2(\Omega)$ and $L^2(\Omega)^3$ with square integrable gradients, divergences, and curls, respectively. We define their subspaces $L_0^2(\Omega)$, $H_0^1(\Omega)$, $\bm{H}_0(\mathrm{div};\Omega)$, $\bm{H}_0(\mathrm{curl};\Omega)$ with vanishing mean values on $\Omega$, vanishing traces, normal traces, and tangential traces on $\partial\Omega$, respectively. Then, let $Q_h\subset{H}^1(\Omega)$, $\bm{X}_h\subset{H}_0^1(\Omega)^3$, and $M_h\subset{L}_0^2(\Omega)\cap{H}^1(\Omega)$ be the $H^1$-conforming finite element spaces, $\bm{N}_h\subset\bm{H}(\mathrm{div};\Omega)$ be the $\bm{H}(\mathrm{div})$-conforming finite element space, and $\bm{C}_h\subset\bm{H}_0(\mathrm{curl};\Omega)$ be the $\bm{H}(\mathrm{curl})$-conforming finite element space. It is assumed that $\nabla\times\bm{C}_h\subset\bm{N}_h$. We use a uniform partition of the time interval with a time step size $\tau$. Then we define
\begin{equation*}
    \delta_\tau\chi^{n+1}=\frac{\gamma_0\chi^{n+1}-\widehat{\chi}}{\tau},
\end{equation*}
where $\chi$ is a generic variable, and
\begin{equation*}
    \widehat{\chi}=
    \begin{cases}
        \chi^n, & k=1, \\
        2\chi^n-\frac{1}{2}\chi^{n-1}, & k=2,
    \end{cases}
    ; \quad \gamma_0=
    \begin{cases}
        1, & k=1, \\
        \tfrac{3}{2}, & k=2.
    \end{cases}
\end{equation*}
We also recall the extrapolation formulas
\begin{equation*}
    \widetilde{\chi}^{n+1}=
    \begin{cases}
        \chi^n, & k=1, \\
        2\chi^n-\chi^{n-1}, & k=2.
    \end{cases}
\end{equation*}
With these notations, the discrete finite element scheme is given as follows:

\textit{Step} 1: Find $\big(\phi_h^{n+1},\mu_h^{n+1}\big)\in{Q}_h\times{Q}_h$, such that for all $(\psi_h,\omega_h)\in{Q}_h\times{Q}_h$,
\begin{equation*}
    \begin{aligned}
        \big(\delta_\tau\phi_h^{n+1},\omega_h\big)+\frac{1}{\mathrm{Pe}}\big(\nabla\mu_h^{n+1},\nabla\omega_h\big)= & \big(\widetilde{\phi}_h^{n+1}\widetilde{\bm{u}}_h^{n+1},\nabla\omega_h\big) , \\
        \big(\mu_h^{n+1},\psi_h\big)-\mathrm{Cn}^2\big(\nabla\phi_h^{n+1},\nabla\psi_h\big)-\big(\phi_h^{n+1},\psi_h\big)=  & \big((\widetilde{\phi}_h^{n+1})^3-2\widetilde{\phi}_h^{n+1},\psi_h\big).
    \end{aligned}
\end{equation*}
Then we evaluate $\rho_h^{n+1}$, $\eta_h^{n+1}$, and $\sigma_h^{n+1}$ by (exemplarily for $\rho_h^{n+1}$)
\begin{equation*}
    \rho_h^{n+1}=\frac{\rho_+-\rho_-}{2\rho_r}\mathscr{C}\big(\phi_h^{n+1}\big)+\frac{\rho_++\rho_-}{2\rho_r},
\end{equation*}
where we define
\begin{equation*}
    \mathscr{C}\big(\phi_h^{n+1}\big)=
    \begin{cases}
        \phi_h^{n+1}, & \text{if}\;\big|\phi_h^{n+1}\big|\le{1}, \\
        \mathrm{sign}(\phi_h^{n+1}), & \text{otherwise}.
    \end{cases}
\end{equation*}

\textit{Step} 2: Find $\bm{u}_h^{n+1}\in\bm{X}_h$, such that for all $\bm{v}_h\in\bm{X}_h$,
\begin{equation*}
    \begin{aligned}
        &\big(\rho_h^{n+1}\delta_\tau\bm{u}_h^{n+1},\bm{v}_h\big)+\frac{2}{\mathrm{Re}}\big(\eta_h^{n+1}\D(\bm{u}_h^{n+1}),\D(\bm{v}_h)\big)=-\big((\rho_h^{n+1}\widetilde{\bm{u}}_h^{n+1}-\rho_d\tfrac{1}{\mathrm{Pe}}\nabla\mu_h^{n+1})\cdot\nabla\widetilde{\bm{u}}_h^{n+1},\bm{v}_h\big) \\
         &+\big(\widetilde{p}_h^{n+1},\nabla\cdot\bm{v}_h\big)+\frac{1}{\mathrm{We}\mathrm{Cn}}\big(\mu_h^{n+1}\nabla\phi_h^{n+1},\bm{v}_h\big)-\frac{\mathrm{N}}{\mathrm{Rm}}\big(\widetilde{\bm{B}}_h^{n+1}\otimes\widetilde{\bm{B}}_h^{n+1},\D(\bm{v}_h)\big)-\frac{1}{\mathrm{Fr}}\big(\rho_h^{n+1}\bm{e}_3,\bm{v}_h\big).
    \end{aligned}
\end{equation*}

\textit{Step} 3: Find $p_h^{n+1}\in{M}_h$, such that for all $q_h\in{M}_h$,
\begin{equation*}
    \big(\nabla({p}_h^{n+1}-{p}_h^n),\nabla{q}_h\big)=-\frac{\gamma_0\vartheta}{\tau}\big(\nabla\cdot\bm{u}_h^{n+1},q_h\big),
\end{equation*}
with $\vartheta=\tfrac{\min(\rho_+,\rho_-)}{\rho_r}$.

\textit{Step} 4: Find $\bm{A}_h^{n+1}\in\bm{C}_h$, such that for all $\bm{Z}_h\in\bm{C}_h$,
\begin{equation*}
    \big(\sigma_h^{n+1}\delta_\tau\bm{A}_h^{n+1},\bm{Z}_h\big)+\frac{1}{\mathrm{Rm}}\big(\nabla\times\bm{A}_h^{n+1},\nabla\times\bm{Z}_h\big)=\big(\sigma_h^{n+1}\bm{u}_h^{n+1}\times\nabla\times\widetilde{\bm{A}}_h^{n+1},\bm{Z}_h\big).
\end{equation*}

\textit{Step} 5: Find $\bm{B}_h^{n+1}\in\bm{N}_h$, such that for all $\bm{D}_h\in\bm{N}_h$,
\begin{equation*}
    \big(\bm{B}_h^{n+1},\bm{D}_h\big)=\big(\nabla\times\bm{A}_h^{n+1},\bm{D}_h\big).
\end{equation*}

\begin{remark}
    In the above scheme, we have exploited a linearization technique in \cite{2015_Shen} and the pressure penalty method developed by Guermond and Salgado \cite{2009_Guermond}. While this scheme offers a straightforward implementation and exactly preserves the divergence-free constraint on the magnetic field $\bm{B}$, it neither satisfies the fluid incompressibility constraint nor admits a discrete stability estimate. Designing an algorithm that preserves these properties simultaneously remains a significant challenge and hence is beyond the scope of this work.
\end{remark}

\begin{remark}
    By its physical interpretation, the order parameter $\phi$ should take values in $[-1,1]$. However, the present formulation cannot theoretically guarantee this bound \cite{2019_Giorgini}. In computations, the cut-off value $\mathscr{C}(\phi_h^{n+1})$ defined above is used when evaluating material properties $\rho_h^{n+1}$, $\eta_h^{n+1}$, and $\sigma_h^{n+1}$. This treatment guarantees that these material properties remain within their physically admissible ranges. It is especially important for simulations with large parameter ratios, and is a common strategy in simulations \cite{2015_Shen, 2024_Wang}.
\end{remark}

\begin{remark}
    In computations, we fix the Froude number $\mathrm{Fr}=1$, which gives $u_r=\sqrt{gL_r}$. The Reynolds number and the Weber number now become the Archimedes number $\mathrm{Ar}=\tfrac{\rho_+\sqrt{gL_r^3}}{\eta_+}=\mathrm{Re}$ and the Eötvös number $\mathrm{Eo}=\tfrac{L_r^2\rho_+g}{\lambda}=\tfrac{2\sqrt{2}}{3}\mathrm{We}$, where we have chosen $\chi_r=\chi_+$ with $\chi=\{\rho,\eta,\sigma\}$. We also employ the Hartmann number $\mathrm{Ha}=L_rB_r\sqrt{\tfrac{\sigma_+}{\eta_+}}=\sqrt{\mathrm{N}\mathrm{Re}}$ and the magnetic Prandtl number $\mathrm{Pm}=\tfrac{\eta_+\mu_0\sigma_+}{\rho_+}=\tfrac{\mathrm{Rm}}{\mathrm{Re}}$ in place of the Stuart number $\mathrm{N}$ and the magnetic Reynolds number $\mathrm{Rm}$. With the scaling law $\tfrac{1}{\mathrm{Pe}}=3\mathrm{Cn}$ \cite{2013_Magaletti}, the system is fully characterized by eight dimensionless quantities: $\tfrac{\rho_+}{\rho_-}$, $\tfrac{\eta_+}{\eta_-}$, $\tfrac{\sigma_+}{\sigma_-}$, $\mathrm{Cn}$, $\mathrm{Ar}$, $\mathrm{Eo}$, $\mathrm{Ha}$, $\mathrm{Pm}$.
\end{remark}

\section{Numerical examples}
\label{section_3}

In this section, we present numerical examples demonstrating the magnetic damping effect on bubble dynamics. In metallurgical processes, gas bubbles (e.g., argon and nitrogen) are frequently injected into liquid metals to stir and homogenize the melt or to prevent nozzle clogging. Meanwhile, magnetic fields provide a non-intrusive means of bubble control \cite{2024_Gou, 2014_Zhang, 2018_Zhang}. Here, we consider a single lighter bubble with an initial diameter $L_r=0.5$ rising in a vertical column of heavier liquid. The computational domain is $(0,2L_r)\times(0,4L_r)$ in two dimensions, and $(0,2L_r)^2\times(0,4L_r)$ in three dimensions, with no-slip and insulating boundary conditions on all walls. The initial center of the bubble is positioned at $(L_r,L_r)$ and $(L_r,L_r,L_r)$ in two and three dimensions, respectively. The initial order parameter profile is prescribed as
\begin{equation*}
    \phi^0=\tanh\frac{d(\x)}{\sqrt{2}\mathrm{Cn}},
\end{equation*}
where $d(\x)$ denotes the signed distance function to the bubble interface, chosen to be positive outside the bubble. The fluid is initially at rest, i.e., $\bm{u}^0=\bm{0}$. In all simulations, the initial magnetic field $\bm{B}^0$ is prescribed as the imposed magnetic field, and the corresponding initial vector potential $\bm{A}^0$ is chosen as specified earlier. In three-dimensional examples, uniform magnetic fields are applied in both the vertical and horizontal directions that orient antiparallel and perpendicular to gravity, respectively. The two-dimensional simulations are restricted to a vertical field antiparallel to the direction of gravity. The magnetic Prandtl number is prescribed as either $\mathrm{Pm}=10^{-2}$ or $\mathrm{Pm}=10$, corresponding to the low and high magnetic Reynolds number regimes. The relevant physical parameters are $\tfrac{\rho_+}{\rho_-}=10^3$, $\tfrac{\eta_+}{\eta_-}=10^2$, $\tfrac{\sigma_+}{\sigma_-}=10^4$, $\mathrm{Ar}=35$, $\mathrm{Eo}=125$. For a quantitative comparison, we use the centroid and rise velocity defined by
\begin{equation*}
    \frac{\int_{\phi<0}x_d\d\x}{\int_{\phi<0}\d\x}
    \qquad \mbox{and} \qquad
    \frac{\int_{\phi<0}u_d\d\x}{\int_{\phi<0}\d\x},
\end{equation*}
where $d=2,3$ and $u_d$ is the $d$th component of $\bm{u}=(u_1,u_2,u_3)$.

Our numerical implementation relies on the \texttt{deal.II} finite element library \cite{2024_Africa}. We apply the lowest-order elements for all finite element spaces, specifically using a tensor-product continuous Lagrange element $\mathbb{Q}_1$ for $\phi_h$, $\mu_h$, $\bm{u}_h$, and $p_h$, the Raviart--Thomas element $\mathbb{RT}_0$ for $\bm{B}_h$, and the N\'{e}d\'{e}lec element $\mathbb{N}_0$ for $\bm{A}_h$. The resulting linear systems are solved using the flexible generalized minimal residual (FGMRES) algorithm, preconditioned by either additive Schwarz methods or algebraic multigrid methods. For brevity, the details of preconditioning are omitted as they are not the primary focus of this study. Additionally, we use an absolute iteration tolerance of $\varepsilon_{iter}=10^{-10}$ in all examples. All numerical results reported below are rescaled to the physical dimensional scale for presentation, and, unless otherwise specified, all simulations are performed up to the physical final time $t=3$.

\subsection{Convergence test}

In this subsection, we investigate the effect of interface thicknesses on numerical results to examine the asymptotic behavior of our diffuse interface model. Five two-dimensional experiments are performed with $\mathrm{Cn}=0.04,0.02,0.01,0.005$, and $0.0025$. To adequately resolve the interfacial dynamics, we employ a fixed locally refined mesh for each Cahn number. Specifically, a uniform background mesh with mesh size $h$ is first generated, followed by one level of local refinement in the subdomain $(0.15,0.85)\times(0.2,1.5)$. As the Cahn number decreases, both the mesh resolution and the time step size $\tau$ are refined accordingly, as summarized in Table \ref{comparison_Cn_0}. The remaining dimensionless parameters are fixed at $\mathrm{Pm}=10^{-2}$ and $\mathrm{Ha}=15$.

\begin{table}[!htbp]
    \centering
    \caption{Spatial and temporal resolutions in the convergence test.}
    \begin{tabular}{llll}
        \toprule
        $\mathrm{Cn}$ & $\tau$ & $h$ & cells \\
        \midrule
        0.04 & $2\times10^{-3}$ & $2^{-5}$ & 4,820 \\
        0.02 & $1\times10^{-3}$ & $2^{-6}$ & 19,148 \\
        0.01 & $5\times10^{-4}$ & $2^{-7}$ & 77,588 \\
        0.005 & $2.5\times10^{-4}$ & $2^{-8}$ & 310,892 \\
        0.0025 & $1.25\times10^{-4}$ & $2^{-9}$ & 1,239,572 \\
        \bottomrule
    \end{tabular}
    \label{comparison_Cn_0}
\end{table}

\begin{figure}[!htbp]
    \centering
    \subfloat[Final bubble shapes]{
        \includegraphics[width=0.3\linewidth]{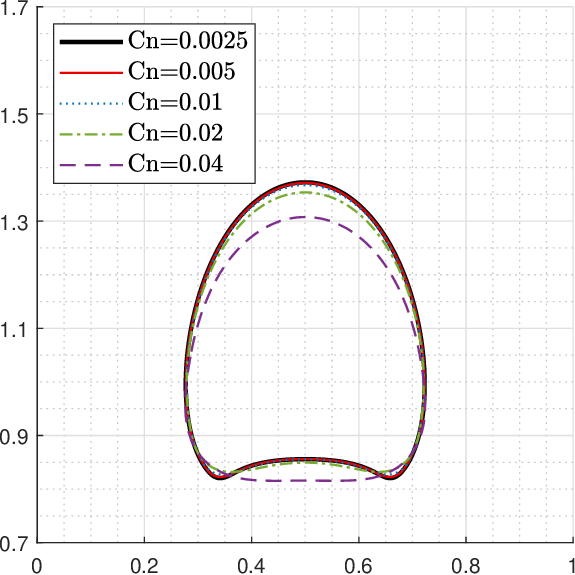}
    }
    \\
    \subfloat[Time evolution of centroid and rise velocity]{
        \includegraphics[width=0.6\linewidth]{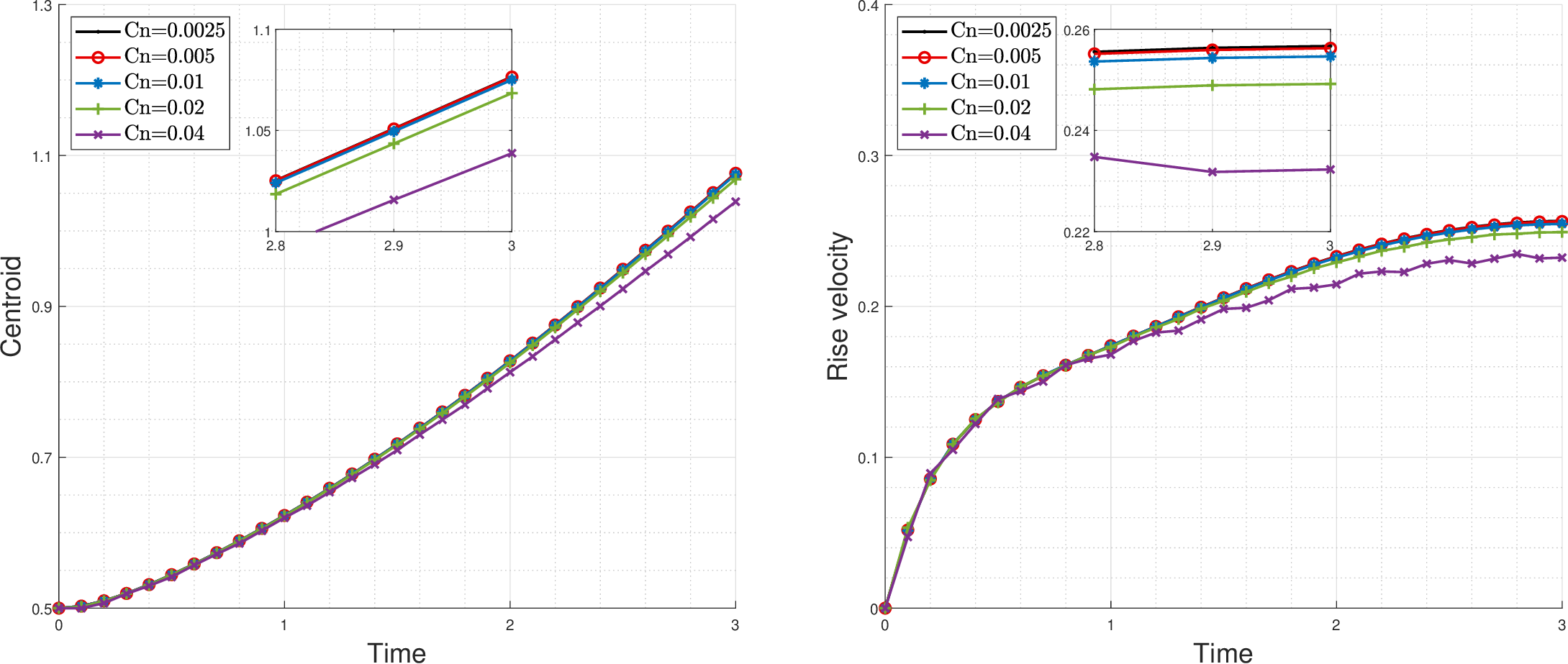}
    }
    \caption{Effect of the Cahn number on the computed bubble dynamics.}
    \label{comparison_Cn_1}
\end{figure}

Figure \ref{comparison_Cn_1} shows the final bubble shapes for different Cahn numbers, together with the time evolution of the corresponding benchmark quantities. It is evident that the numerical results converge to a limiting solution as the Cahn number decreases. For a more quantitative assessment, Table \ref{comparison_Cn_2} summarizes the relative differences and estimated convergence orders of the final centroid and rise velocity, using the result obtained with $\mathrm{Cn}=0.0025$ as the reference solution. The results indicate an approximately second-order convergence for both benchmark quantities. In particular, no pronounced difference in either the bubble shape or the benchmark quantities is observed between the cases $\mathrm{Cn}=0.005$ and $\mathrm{Cn}=0.0025$.

\begin{table}[!htbp]
    \centering
    \caption{Final benchmark quantities and estimated convergence order.}
    \begin{tabular}{llcclcc}
        \toprule
        \multirow{2}{*}{$\mathrm{Cn}$} & \multicolumn{3}{c}{Centroid} & \multicolumn{3}{c}{Rise velocity} \\ \cmidrule(lr){2-4} \cmidrule(lr){5-7}
        & value & error & order & value & error & order \\
        \midrule
        0.04 & 1.0387 & $3.52\times{10}^{-2}$ & & 0.2323 & $9.50\times{10}^{-2}$ & \\
        0.02 & 1.0684 & $7.56\times{10}^{-3}$ & 2.22 & 0.2492 & $2.92\times{10}^{-2}$ & 1.70 \\
        0.01 & 1.0750 & $1.45\times{10}^{-3}$ & 2.38 & 0.2546 & $8.03\times{10}^{-3}$ & 1.86 \\
        0.005 & 1.0764 & $1.99\times{10}^{-4}$ & 2.87 & 0.2563 & $1.69\times{10}^{-3}$ & 2.25 \\
        0.0025 & 1.0766 & & & 0.2567 & & \\
        \bottomrule
    \end{tabular}
    \label{comparison_Cn_2}
\end{table}

In the subsequent two-dimensional simulations, we employ $\mathrm{Cn}=0.005$ together with a background mesh size $h=2^{-7}$ and a locally refined region $(0.15,0.85)\times(0.2,1.5)$ to reduce the computational cost. For the three-dimensional simulations, we use $\mathrm{Cn}=0.01$ together with a background mesh size $h=2^{-6}$ and a locally refined region $(0.15,0.85)^2\times(0.2,1.8)$, which gives 1,906,592 cells. The time step size is fixed at $\tau=1\times{10}^{-3}$ in all subsequent simulations.

\subsection{Comparison with inductionless MHD and sharp interface ALE computations}

In most laboratory experiments and industrial processes, the magnetic Reynolds number $\mathrm{Rm}$ is typically low \cite{2017_Davidson}, and the inductionless MHD model is therefore widely used as an efficient formulation for the simulation of resistive MHD free surface flows \cite{2019_Herreman, 2004_Morley, 2018_Pan, 2014_Zhang, 2018_Zhang}. It is therefore of interest to examine whether, in the low $\mathrm{Rm}$ regime, the present resistive MHD model produces results consistent with those of the inductionless MHD model. In the computation, the evolving interface and fluid unknowns are computed by the same diffuse interface formulation developed in this work. Meanwhile, for the electromagnetic fields, we follow the procedure in \cite{2014_Zhang}. Specifically, the electric potential is obtained from the electric potential Poisson equation, and the current density is subsequently projected into the $\bm{H}(\mathrm{div})$-conforming finite element space using Ohm's law and charge-conservation constraint, thereby preserving a divergence-free current density at the discrete level.

On the other hand, although the preceding subsection demonstrates numerically that the diffuse interface solutions approach a stable limiting behavior as the Cahn number decreases, a numerical verification is still necessary to validate whether this limiting behavior is consistent with the corresponding sharp interface dynamics. Such a reference solution is obtained using an ALE approach, in which the moving sharp interface is explicitly represented by a body-fitted mesh. The evolving interface and fluid unknowns are computed using the Barrett--Garcke--N\"{u}rnberg type scheme employed in \cite{2026_Wang}, whereas the electromagnetic fields are discretized using the potential formulation developed in the present work.

\begin{figure}[!htbp]
    \centering
    \subfloat[Final bubble shapes]{
        \includegraphics[width=0.3\linewidth]{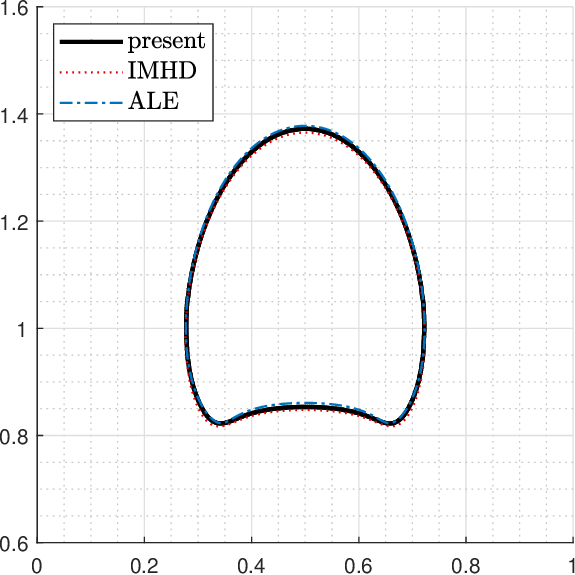}
    }
    \\
    \subfloat[Time evolution of centroid and rise velocity]{
        \includegraphics[width=0.6\linewidth]{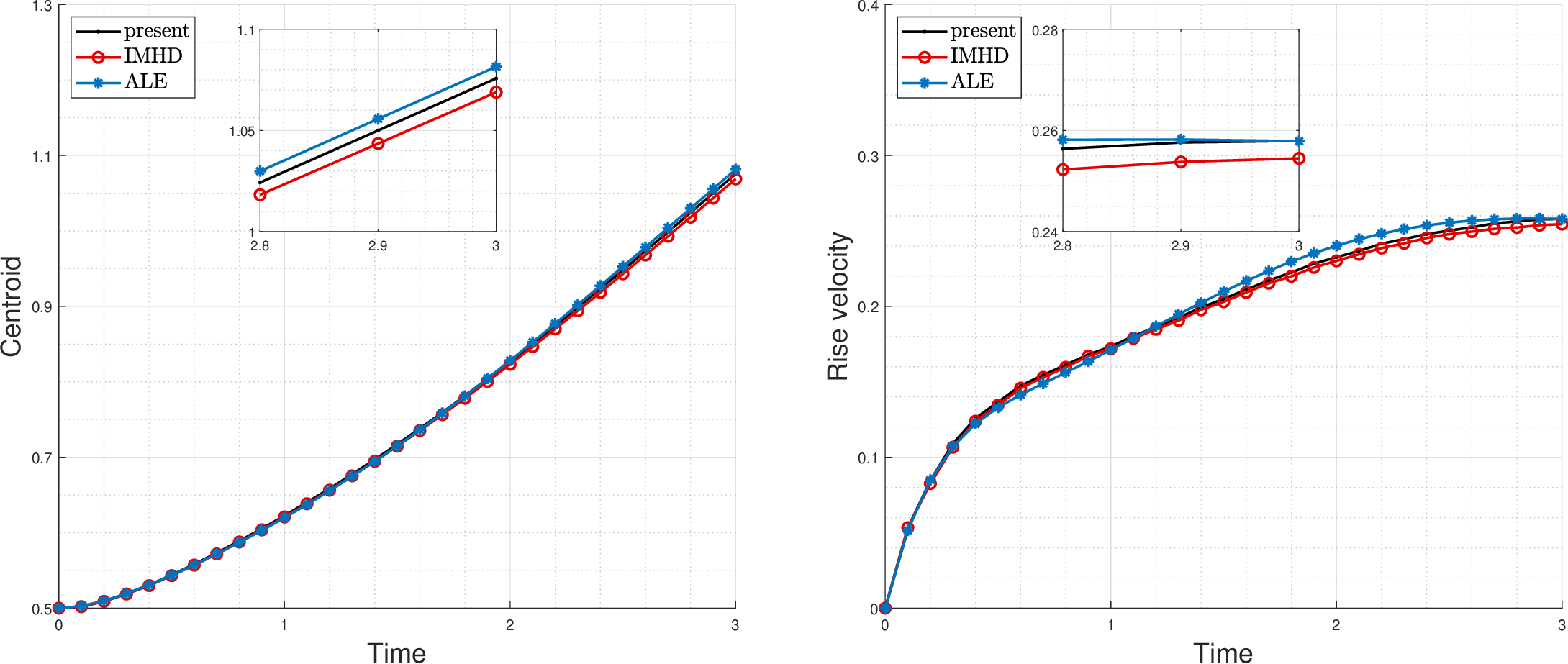}
    }
    \caption{Bubble dynamics obtained from the present, inductionless MHD, and ALE computations. In the legends, present, IMHD, and ALE denote these three results, respectively.}
    \label{comparison_method}
\end{figure}

In this comparison, the relevant dimensionless parameters are $\mathrm{Pm}=10^{-2}$ and $\mathrm{Ha}=15$. Figure \ref{comparison_method} depicts the final bubble shape and the associated benchmark quantities. Close agreement is observed between the present resistive MHD computation and the inductionless MHD computation, as well as between the diffuse interface computation and the sharp interface ALE computation. The former comparison verifies the consistency of the present electromagnetic formulation with the inductionless approximation in the low $\mathrm{Rm}$ regime, whereas the latter validates the diffuse interface approximation against the corresponding sharp interface dynamics. These results indicate the accuracy and robustness of our method in simulating MHD interfacial phenomena.

\subsection{Results in low magnetic Reynolds number regime}

\begin{figure}[!htbp]
    \centering
    \subfloat[$\mathrm{Ha}=5$]{
        \includegraphics[width=0.2\linewidth]{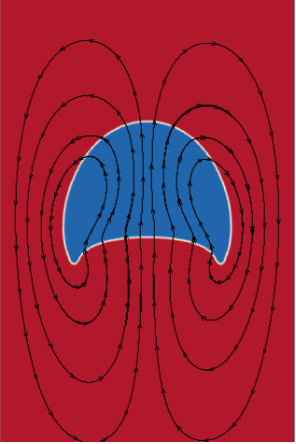}
    }
    \hspace{30pt}
    \subfloat[$\mathrm{Ha}=15$]{
        \includegraphics[width=0.2\linewidth]{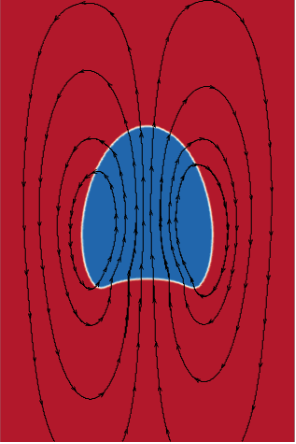}
    }
    \hspace{30pt}
    \subfloat[$\mathrm{Ha}=25$]{
        \includegraphics[width=0.2\linewidth]{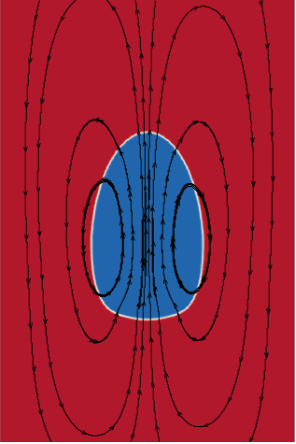}
    }
    \caption{Final order parameter profiles and streamlines in two-dimensional tests at $\mathrm{Pm}=10^{-2}$.}
    \label{low_Rm_2D_shape}
\end{figure}

\begin{figure}[!htbp]
    \centering
    \includegraphics[width=0.8\linewidth]{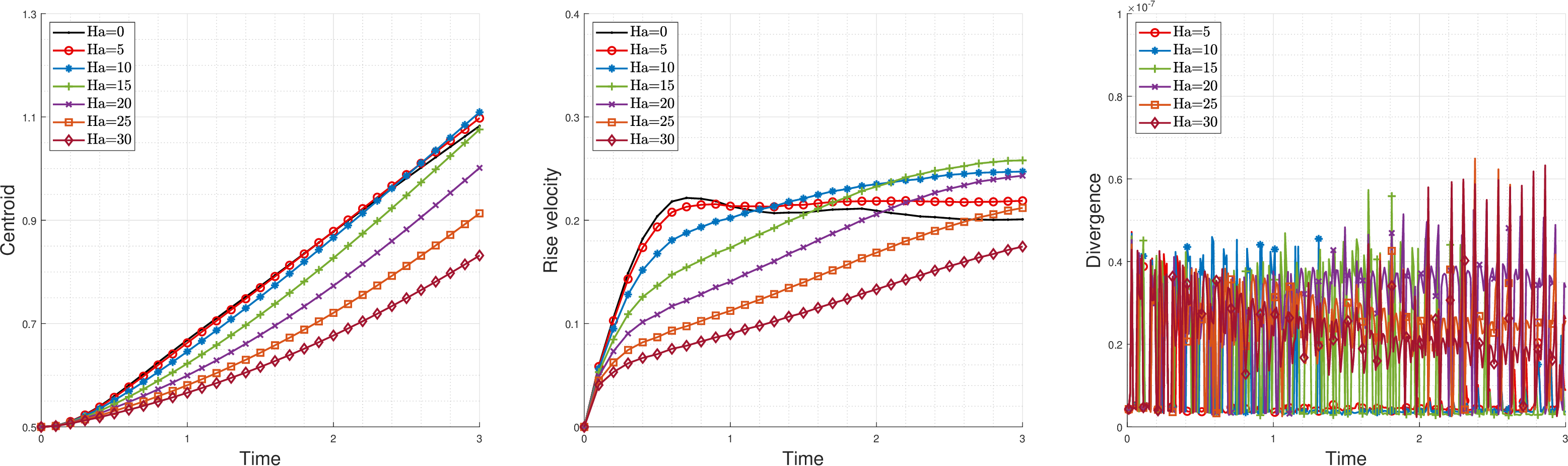}
    \caption{Centroid, rise velocity, and divergence of the magnetic field in two-dimensional tests at $\mathrm{Pm}=10^{-2}$.}
    \label{low_Rm_2D_quantity}
\end{figure}

This subsection presents numerical results for $\mathrm{Pm}=10^{-2}$. We first perform two-dimensional simulations with different magnetic strengths. A comparison of the final order parameter profile under different Hartmann numbers is presented in Figure \ref{low_Rm_2D_shape}. Notably, the bubble undergoes elongation along the vertical direction, giving a more pronounced elliptical profile as the magnetic field intensifies. Meanwhile, the swirling vortices become more regular, with streamlines aligning parallel to the magnetic field lines. For a detailed quantitative assessment, Figure \ref{low_Rm_2D_quantity} illustrates the time evolution of the centroid and rise velocity. Remarkably, the rise velocity depends non-monotonically on the magnetic field strength, and a moderate field is found to enhance the rise velocity, leading to a higher centroid position.

\begin{figure}[!htbp]
    \centering
    \subfloat[Vertical magnetic fields]{
        \includegraphics[width=0.4\linewidth]{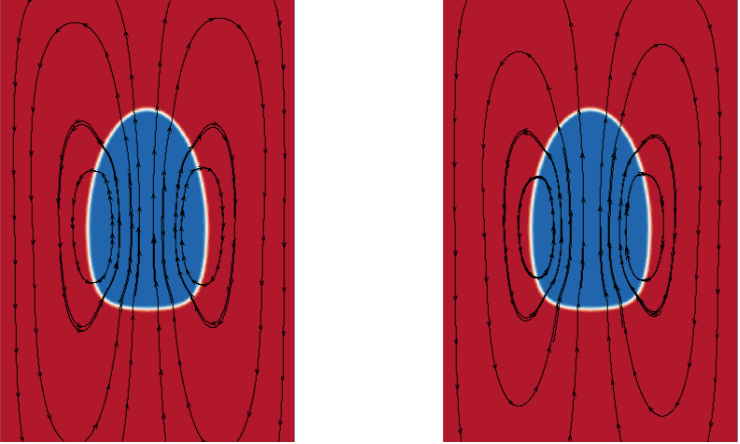}
    }
    \hspace{30pt}
    \subfloat[Horizontal magnetic field]{
        \includegraphics[width=0.4\linewidth]{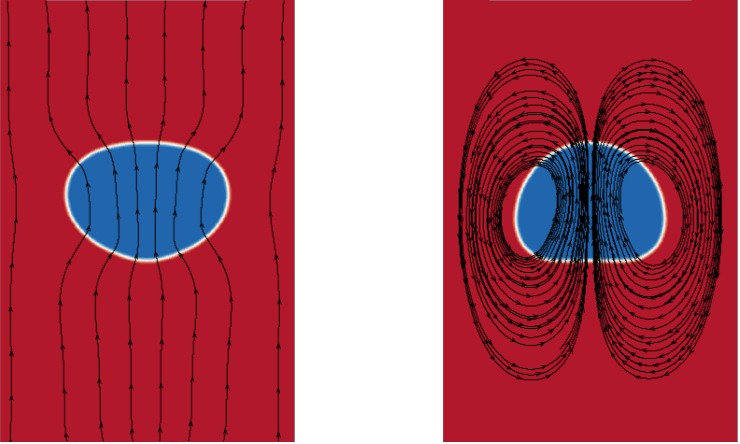}
    }
    \caption{Final order parameter profiles and streamlines for three-dimensional tests at $\mathrm{Ha}=30$ and $\mathrm{Pm}=10^{-2}$. Left and right panels show cross sections parallel and perpendicular to the horizontal magnetic field direction, respectively.}
    \label{low_Rm_3D_shape}
\end{figure}

\begin{figure}[!htbp]
    \centering
    \includegraphics[width=0.8\linewidth]{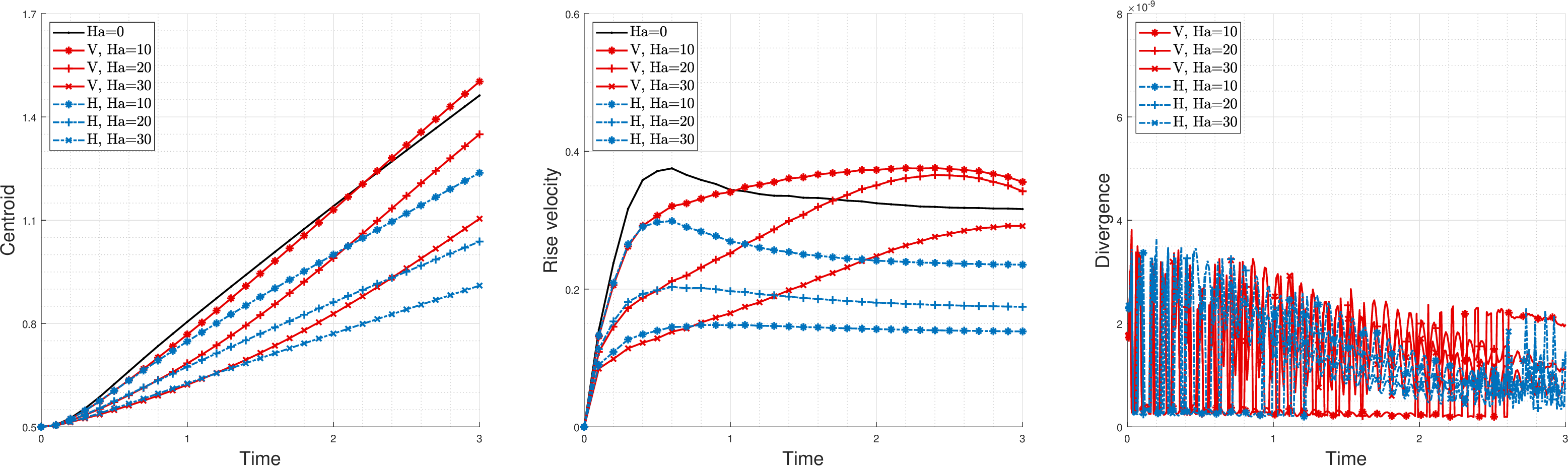}
    \caption{Centroid, rise velocity, and divergence of the magnetic field in three-dimensional tests at $\mathrm{Pm}=10^{-2}$. In the legends, ``V'' and ``H'' denote vertical and horizontal applied magnetic fields, respectively.}
    \label{low_Rm_3D_quantity}
\end{figure}

\begin{figure}[!htbp]
	\centering
	\subfloat[Vertical magnetic fields]{
		\includegraphics[width=0.8\linewidth]{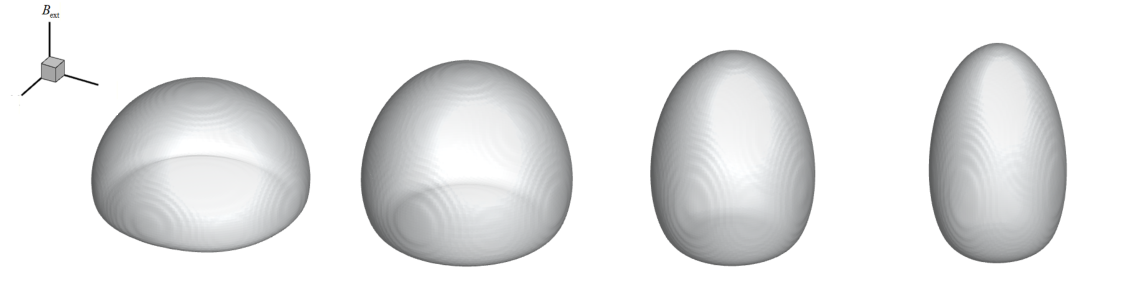}
	}
    \\
	\subfloat[Horizontal magnetic fields]{
		\includegraphics[width=0.8\linewidth]{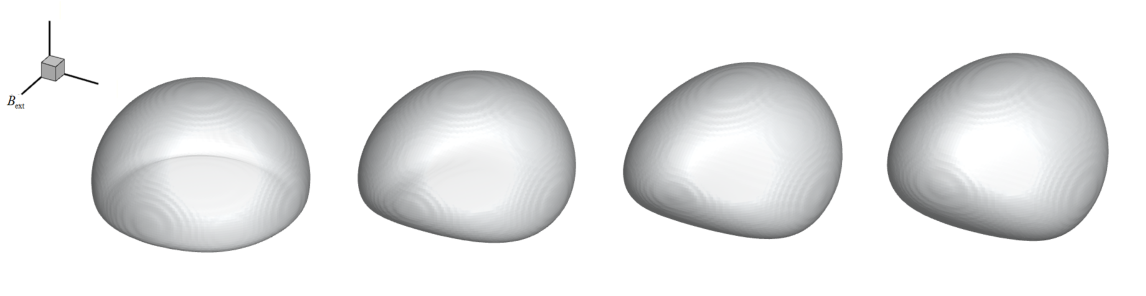}
	}
	\caption{Final three-dimensional bubble shapes for $\mathrm{Pm}=10^{-2}$ with $\mathrm{Ha}=0,10,20$, and $30$. In the legends, $B_\text{ext}$ denotes the direction of the external imposed magnetic field.}
	\label{3D_shape}
\end{figure}

We now turn our attention to three-dimensional simulations. Figure \ref{low_Rm_3D_shape} shows cross sections of final order parameter profiles and streamlines at $\mathrm{Ha}=30$, while the corresponding benchmark quantities are depicted in Figure \ref{low_Rm_3D_quantity}. Under a vertical magnetic field, the observed flow dynamics mirror those in the two-dimensional examples. This can be attributed to the fact that only the velocity components perpendicular to the magnetic field (i.e., the horizontal components) contribute to the induced current. The resulting Lorentz force primarily acts to suppress horizontal motion, leading to an axisymmetric deformation of the bubble. In contrast, a horizontal magnetic field induces a pronounced anisotropic flow structure. In this case, the induced current is generated by velocity components normal to the field direction, resulting in a spatially non-uniform Lorentz force distribution. As a consequence, vortical structures perpendicular to the magnetic field are strongly damped, and the flow becomes highly directional; detailed three-dimensional bubble shapes are shown in Figure \ref{3D_shape}. In particular, the streamlines in the plane normal to the magnetic field tend to align opposite to gravity. Furthermore, both the centroid and rise velocity exhibit a monotonic decrease with increasing magnetic intensity. These numerical observations are in qualitative agreement with previous findings reported in \cite{2014_Zhang}.

In addition, the divergence of $\bm{B}_h$ is maintained at $\mathcal{O}\big(10^{-8}\big)$ in two dimensions and $\mathcal{O}\big(10^{-9}\big)$ in three dimensions, illustrating the robust conservation property of our algorithm in the low $\mathrm{Rm}$ regime; notice that we adopt a tolerance of $\varepsilon_{iter}=10^{-10}$ for FGMRES iterations.

\subsection{Results in high magnetic Reynolds number regime}

\begin{figure}[!htbp]
    \centering
    \subfloat[$\mathrm{Ha}=25$]{
        \includegraphics[width=0.2\linewidth]{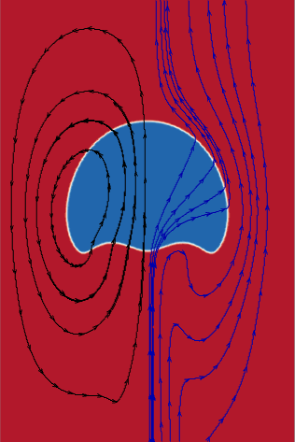}
    }
    \hspace{30pt}
    \subfloat[$\mathrm{Ha}=30$]{
        \includegraphics[width=0.2\linewidth]{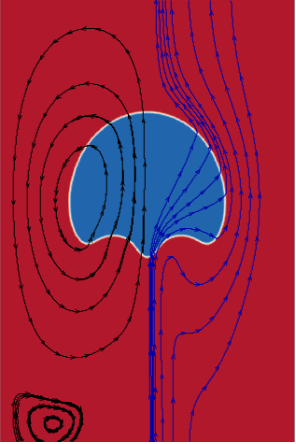}
    }
    \hspace{30pt}
    \subfloat[$\mathrm{Ha}=35$]{
        \includegraphics[width=0.2\linewidth]{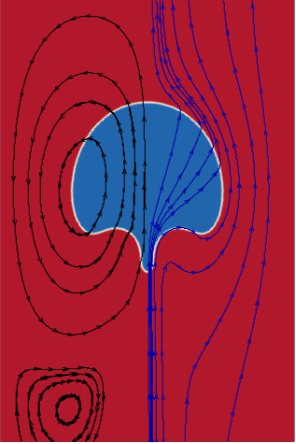}
    }
    \caption{Final order parameter profiles in two-dimensional tests at $\mathrm{Pm}=10$, with streamlines in black and magnetic field lines in red.}
    \label{high_Rm_2D_shape}
\end{figure}

\begin{figure}[!htbp]
    \centering
    \includegraphics[width=0.8\linewidth]{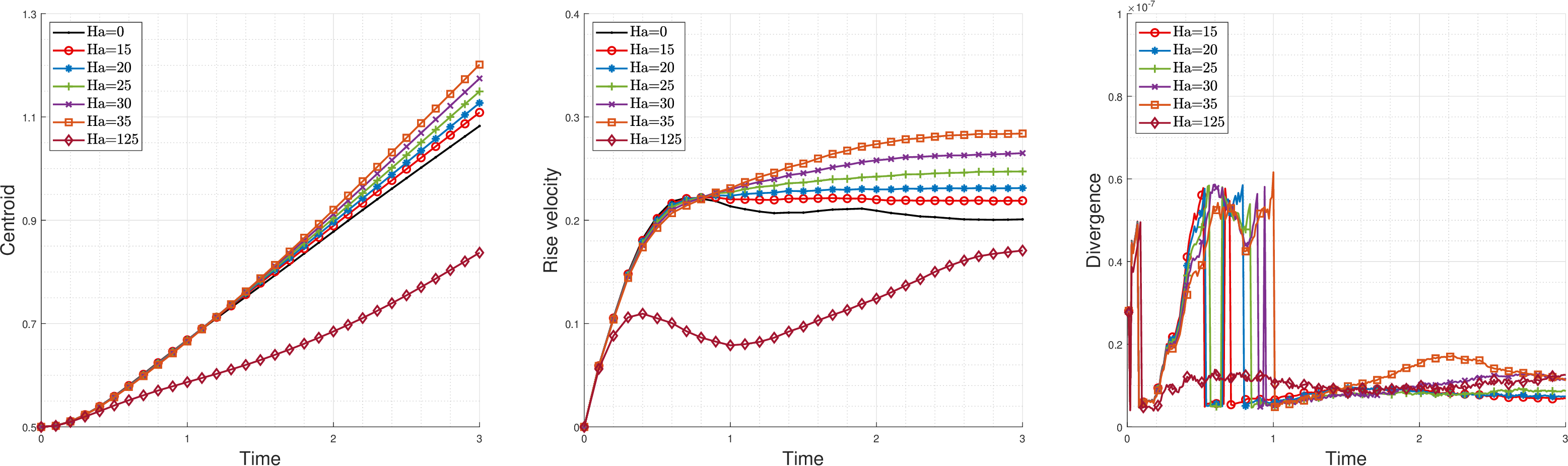}
    \caption{Centroid, rise velocity, and divergence of the magnetic field in two-dimensional tests at $\mathrm{Pm}=10$.}
    \label{high_Rm_2D_quantity}
\end{figure}

This subsection presents numerical results for $\mathrm{Pm}=10$. Given that a vertical magnetic field induces an axisymmetric fluid motion, we restrict our investigation to two-dimensional simulations for this field. Figure \ref{high_Rm_2D_shape} illustrates the final order parameter profiles under various magnetic strengths, with streamlines and magnetic field lines superimposed. Although the vortices become more regular as the magnetic field strength increases, the magnetic damping effect is less pronounced. Meanwhile, the non-monotonic dependence of the rise velocity on the magnetic-field strength is weaker, so that an intense field is required to inhibit the bubble motion, as detailed in Figure \ref{high_Rm_2D_quantity}. Moreover, the magnetic flux freezing effect becomes increasingly significant, causing the magnetic field lines to become nearly frozen into the flow both within and surrounding the bubble. This induces the formation of a protrusion at the bubble's base that remains absent in the low-$\mathrm{Rm}$ regime. Consequently, secondary vortices are generated beneath the bubble.

\begin{figure}[!htbp]
    \centering
    \subfloat[Magnetic field lines]{
        \includegraphics[width=0.4\linewidth]{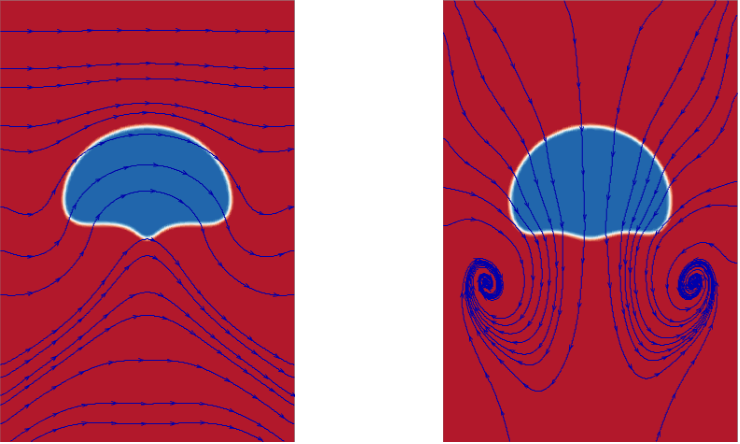}
    }
    \hspace{30pt}
    \subfloat[Streamlines]{
        \includegraphics[width=0.4\linewidth]{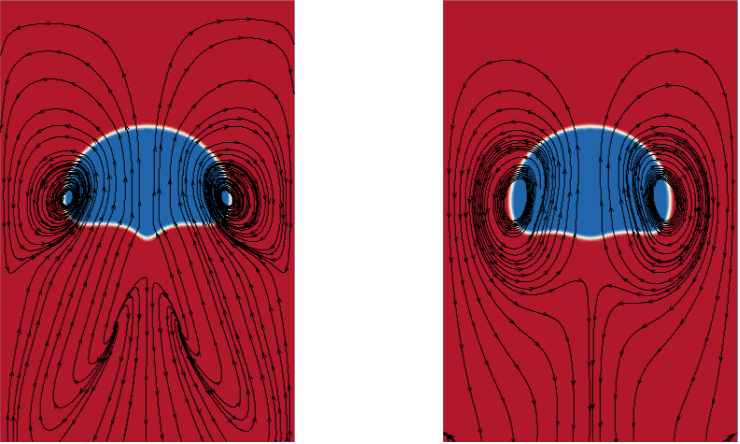}
    }
    \caption{Three-dimensional magnetic field lines and streamlines for a horizontal magnetic field at $\mathrm{Ha}=30$ and $\mathrm{Pm}=10$. Left and right panels show cross sections parallel and perpendicular to the horizontal magnetic field direction, respectively.}
    \label{high_Rm_3D_shape}
\end{figure}

\begin{figure}[!htbp]
    \centering
    \includegraphics[width=0.8\linewidth]{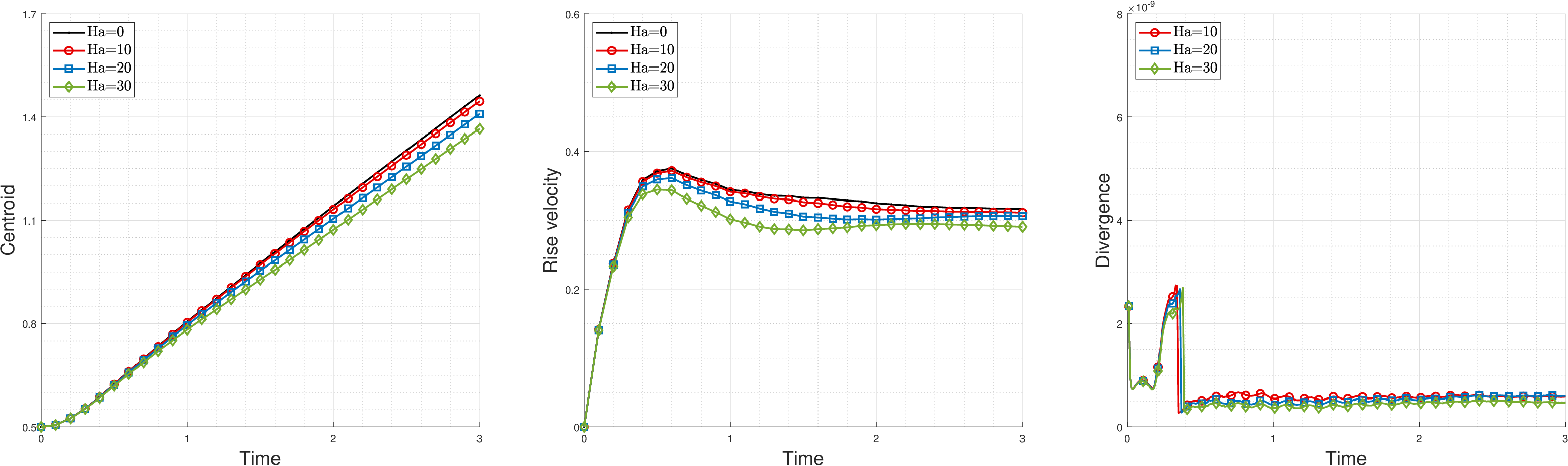}
    \caption{Centroid, rise velocity, and divergence of the magnetic field in a three-dimensional test with a horizontal magnetic field at $\mathrm{Pm}=10$.}
    \label{high_Rm_3D_quantity}
\end{figure}

Next, we focus on the dynamics under a horizontal magnetic field. Figure \ref{high_Rm_3D_shape} presents the magnetic field lines and streamlines at $\mathrm{Ha}=30$. The strong coupling between the magnetic field and the fluid flow leads to highly perturbed magnetic field lines. While the protrusion at the bubble's base remains observable, the damping effect does not yield a distinctly anisotropic bubble shape as observed in the low $\mathrm{Rm}$ regime; see also Figure \ref{high_Rm_3D_quantity} for a quantitative measurement. Furthermore, the flow field becomes increasingly complex, with secondary vortices also emerging beneath the bubble in the parallel plane.

Also, the divergence of $\bm{B}_h$ is maintained at $\mathcal{O}\big(10^{-8}\big)$ in two dimensions and $\mathcal{O}\big(10^{-9}\big)$ in three dimensions, demonstrating the robust conservation property of our algorithm in high $\mathrm{Rm}$ regime.

\section{Concluding remarks}
\label{section_4}

In this paper, we propose a diffuse interface model to simulate the incompressible resistive MHD free surface flows. The formal convergence of the diffuse interface system to the sharp interface model is established by the method of matched asymptotic expansions. For the purpose of obtaining numerical solutions, we develop a fully decoupled linear finite element algorithm that preserves the divergence-free constraint on the magnetic field. The magnetic damping effects are investigated through rising bubble dynamics. Numerical tests demonstrate the accuracy and efficiency of our approach.

Although the proposed model and developed algorithm are capable of effectively capturing such MHD phenomena as the magnetic damping effect in free surface flows, our investigation on this topic remains relatively limited. We identify several directions for future research. First, it is essential to take the thermocapillary effect into account, in view of practical industrial applications \cite{2017_Davidson, 2014_Zhang, 2018_Zhang}. Another direction for future investigation lies in the investigation of other MHD phenomena using the developed method, such as metal pad rolling in stably stratified liquid layers \cite{2003_Gerbeau, 2019_Herreman, 2018_Tucs}. Finally, it is highly desirable to develop structure-preserving numerical algorithms that preserve fundamental physical properties for the MHD free surface flow problems.

\section*{Acknowledgments}

The authors would like to thank Prof. Wenyu Lei at UESTC for the fruitful discussions on numerical implementations using \texttt{deal.II}. The research is supported in part by the National Natural Science Foundation of China (Nos. 12271082, 12471371, 12431015, 62231016), the Foundation for Innovative Research Groups of the Natural Science Foundation of Sichuan Province (No. 2025NSFTD0004), and the Fundamental Research Funds for the Central Universities (No. ZYGX2025XJ029).





\appendix
\section{Details of the numerical model}

\subsection{Derivation of the diffuse interface model}
\label{derivation}

We present here the derivation of system \eqref{phase} by means of Onsager's variational principle and the laws of thermodynamics. First, we assume partial mixing of the two fluids in a narrow interfacial zone $\Gamma_\epsilon(t)$ of thickness $\epsilon$ under the volume incompressibility assumption. We need to define the unknowns in $\Gamma_\epsilon(t)$ and develop the corresponding governing equations. Hence, let $W(t)\subset\Omega$ be an arbitrary control volume, and let $M$ and $\rho$ denote its mass and density, respectively. Then, denoting the masses and densities of the two fluids in $W(t)$ by $M_\pm$ and $\varrho_\pm$, it holds that
 \begin{equation*}
    \varrho_\pm=\frac{M_\pm}{|W(t)|}\in[0,\rho_\pm]
    \qquad \mbox{and} \qquad
    \rho=\frac{M}{|W(t)|}=\frac{M_++M_-}{|W(t)|}=\varrho_++\varrho_-.
\end{equation*}
In turn, the volume incompressibility assumption leads to
\begin{equation*}
    1=\frac{\frac{M_+}{\rho_+}+\frac{M_-}{\rho_-}}{|W(t)|}=\frac{\varrho_+}{\rho_+}+\frac{\varrho_-}{\rho_-}.
\end{equation*}
We define the order parameter for the mixture as
\begin{equation*}
    \phi(\x,t)=\frac{\varrho_+}{\rho_+}-\frac{\varrho_-}{\rho_-}=
    \begin{cases}
        1,  & \x\in\Omega_+(t),\\
        -1, & \x\in\Omega_-(t),
    \end{cases}
\end{equation*}
which indicates that
\begin{equation*}
    \rho=\frac{\rho_+-\rho_-}{2}\phi+\frac{\rho_++\rho_-}{2}.
\end{equation*}
Apparently, $|\phi(\x,t)|<1$ in $\Gamma_\epsilon(t)$. Supposing that the fluids move with different velocities $\bm{u}_\pm$, the continuity equations read
\begin{equation*}
    \partial_t\varrho_\pm+\nabla\cdot(\varrho_\pm\bm{u}_\pm)=0.
\end{equation*}
Then, we introduce the volume-averaged velocity and relative mass flux as
\begin{equation*}
    \bm{u}=\frac{\varrho_+}{\rho_+}\bm{u}_++\frac{\varrho_-}{\rho_-}\bm{u}_-
    \qquad \mbox{and} \qquad
    \bm{j}_\pm=\varrho_\pm(\bm{u}_\pm-\bm{u}).
\end{equation*}
As a consequence, we obtain
\begin{equation*}
    \nabla\cdot\bm{u}=\nabla\cdot\bigg(\frac{\varrho_+}{\rho_+}\bm{u}_++\frac{\varrho_-}{\rho_-}\bm{u}_-\bigg)=-\partial_t\bigg(\frac{\varrho_+}{\rho_+}+\frac{\varrho_-}{\rho_-}\bigg)=0,
\end{equation*}
and
\begin{equation*}
    \partial_t\phi=-\nabla\cdot\bigg(\frac{\varrho_+}{\rho_+}\bm{u}+\frac{\bm{j}_+}{\rho_+}-\frac{\varrho_-}{\rho_-}\bm{u}-\frac{\bm{j}_-}{\rho_-}\bigg)=-\nabla\cdot(\phi\bm{u}+\bm{j}_\phi),
\end{equation*}
where $\bm{j}_\phi=\tfrac{\bm{j}_+}{\rho_+}-\tfrac{\bm{j}_-}{\rho_-}$. Correspondingly, the incompressibility constraint is retained and this implies that mass diffusion is permitted in this method. Also, this gives a non-conservative gravitational force. Hence, the momentum $\rho\bm{u}$ should be transported by $\rho\bm{u}+\rho_d\bm{j}_\phi$, giving
\begin{equation*}
    \rho\partial_t\bm{u}+(\rho\bm{u}+\rho_d\bm{j}_\phi)\cdot\nabla\bm{u}=\nabla\cdot\S-\nabla{p}+\bm{F},
\end{equation*}
where $\S$ is a symmetric second-order tensor, $p$ is the pressure serving as a Lagrange multiplier for the incompressibility constraint, and $\bm{F}$ is a force density. The definition of the electromagnetic fields is independent of the control volume, so we directly invoke the pre-Maxwell equations:
\begin{equation*}
    \begin{aligned}
        \partial_t\bm{B}+\nabla\times\bm{E} & =\bm{0}, \\
        \nabla\times\bm{B} & =\mu_0\bm{J}, \\
        \bm{J} & =\sigma(\phi)(\bm{E}+\bm{u}\times\bm{B}), \\
        \nabla\cdot\bm{B} & =0,\quad\nabla\cdot\bm{J}=0,
    \end{aligned}
\end{equation*}
where the definition of electrical conductivity $\sigma(\phi)$ will be specified later.

Following \cite{2012_Abels, 2006_Qian}, we employ Onsager's variational principle to determine $\bm{j}_\phi$, $\S$, and $\bm{F}$. By introducing the Rayleighian functional
\begin{equation*}
    \mathcal{R}=\frac{\d\mathcal{E}}{\d{t}}+\mathcal{P},
\end{equation*}
where $\mathcal{E}$ and $\mathcal{P}$ are the free energy and dissipation functional, Onsager's variational principle states that the kinetic equation is obtained by minimizing $\mathcal{R}$ \cite{2006_Qian}. Then, we define
\begin{equation*}
    \mathcal{E}=\int_\Omega\bigg(\frac{1}{2}\rho|\bm{u}|^2+\frac{1}{2}\mu_0^{-1}|\bm{B}|^2+f(\phi,\nabla\phi)\bigg)\d\x,
\end{equation*}
where $f(\phi,\nabla\phi)$ is a Helmholtz free energy density. Hence, it holds that
\begin{equation*}
    \frac{\d\mathcal{E}}{\d{t}}=\int_\Omega\big(-\S:\D(\bm{u})+\bm{F}\cdot\bm{u}-\sigma^{-1}(\phi)|\bm{J}|^2-\bm{J}\times\bm{B}\cdot\bm{u}-\mu\bm{u}\cdot\nabla\phi+\bm{j}_\phi\cdot\nabla\mu\big)\d\x,
\end{equation*}
where $\mu=\tfrac{\partial{f}}{\partial\phi}-\nabla\cdot\big(\tfrac{\partial{f}}{\partial\nabla\phi}\big)$ denotes the chemical potential, and here we assume $\tfrac{\partial{f}}{\partial\nabla\phi}\cdot\bm{n}=0$, $\bm{n}\times\bm{E}=\bm{0}$, and $\bm{u}=\bm{0}$ on $\partial\Omega$. Since the first law of thermodynamics reads
\begin{equation*}
    \frac{\d\mathcal{E}}{\d{t}}=-\mathcal{T}\frac{\d\mathcal{S}}{\d{t}}+\frac{\d\mathcal{W}}{\d{t}},
\end{equation*}
where $\mathcal{T}$, $\mathcal{S}$, and $\mathcal{W}$ represent the temperature, entropy, and mechanical work, we are able to identify
\begin{equation*}
    \frac{\d\mathcal{W}}{\d{t}}=\int_\Omega(\bm{F}\cdot\bm{u}-\bm{J}\times\bm{B}\cdot\bm{u}-\mu\bm{u}\cdot\nabla\phi)\d\x.
\end{equation*}
Incorporating gravitational forces, we can write
\begin{equation*}
    \frac{\d\mathcal{W}}{\d{t}}=\int_\Omega\rho\bm{g}\cdot\bm{u}\d\x,
\end{equation*}
which gives
\begin{equation*}
    \bm{F}=\mu\nabla\phi+\bm{J}\times\bm{B}+\rho\bm{g}.
\end{equation*}
To further specify $\S$ and $\bm{j}_\phi$, we introduce the dissipation functional
\begin{equation*}
    \mathcal{P}=\int_\Omega\bigg(\frac{|\S|^2}{4\eta(\phi)}+\frac{|\bm{j}_\phi|^2}{2m_\epsilon(\phi)}+\frac{|\bm{J}|^2}{2\sigma(\phi)}\bigg)\d\x,
\end{equation*}
which accounts for the viscous dissipation, the diffusive dissipation, and the ohmic dissipation, respectively. Hence, we obtain a Rayleighian functional
\begin{equation*}
    \mathcal{R}=\int_\Omega\bigg(-\S:\D(\bm{u})+\rho\bm{g}\cdot\bm{u}+\bm{j}_\phi\cdot\nabla\mu+\frac{|\S|^2}{4\eta(\phi)}+\frac{|\bm{j}_\phi|^2}{2m_\epsilon(\phi)}-\frac{|\bm{J}|^2}{2\sigma(\phi)}\bigg)\d\x.
\end{equation*}
In accordance with Onsager's variational principle, we find by minimizing $\mathcal{R}$ with respect to $\S$ and $\bm{j}_\phi$,
\begin{equation*}
    \S=2\eta\D(\bm{u})
    \qquad \mbox{and} \qquad
    \bm{j}_\phi=-m_\epsilon\nabla{\mu}.
\end{equation*}
In this work, we consider the Helmholtz free energy density in the Ginzburg--Landau form
\begin{equation*}
    f(\phi,\nabla\phi)=\frac{\widehat{\lambda}\epsilon}{2}|\nabla\phi|^2+\frac{\widehat{\lambda}\epsilon^{-1}}{4}\left(\phi^2-1\right)^2.
\end{equation*}
Moreover, we assume that the material properties vary linearly in the mixture, i.e.,
\begin{equation*}
    \eta=\frac{\eta_+-\eta_-}{2}\phi+\frac{\eta_++\eta_-}{2}
    \qquad \mbox{and} \qquad
    \sigma=\frac{\sigma_+-\sigma_-}{2}\phi+\frac{\sigma_++\sigma_-}{2},
\end{equation*}
and hence the model is compatible with the resistive MHD equations in the single fluid regime. This finally leads to system \eqref{phase}.

\subsection{Energy balance law for the diffuse interface model}
\label{law}

We present here the derivation of the energy law for system \eqref{phase}. In the following derivation, all boundary contributions vanish owing to the prescribed boundary conditions. Since $\tfrac{\partial\rho}{\partial\phi}=\rho_d$, the Cahn--Hilliard equation implies
\begin{equation*}
    \partial_t\rho+\nabla\cdot(\rho\bm{u})=\rho_d\nabla\cdot(m_\epsilon\nabla\mu).
\end{equation*}
Accordingly, the momentum equation can be written equivalently in the energy compatible form
\begin{equation*}
    \rho\partial_t\bm{u}+\big(\rho\bm{u}-\rho_dm_\epsilon\nabla\mu\big)\cdot\nabla\bm{u}+\frac{1}{2}(\partial_t\rho+\nabla\cdot(\rho\bm{u})-\rho_d\nabla\cdot(m_\epsilon\nabla\mu))\bm{u}=\nabla\cdot\big(2\eta\D(\bm{u})-p\I\big)+\mu\nabla\phi+\bm{J}\times\bm{B}+\rho\bm{g}.
\end{equation*}
Using integration by parts, we have
\begin{equation*}
    \int_\Omega\big(\rho\bm{u}-\rho_dm_\epsilon\nabla\mu\big)\cdot\nabla\bm{u}\cdot\bm{u}\d\x=-\int_\Omega\frac{1}{2}(\nabla\cdot(\rho\bm{u})-\rho_d\nabla\cdot(m_\epsilon\nabla\mu))|\bm{u}|^2\d\x.
\end{equation*}
Taking the $L^2$ inner product of the momentum equation with $\bm{u}$ therefore yields
\begin{equation*}
    \begin{aligned}
        \frac{\d}{\d{t}}\int_\Omega\frac{1}{2}\rho|\bm{u}|^2\d\x & =\int_\Omega\bigg(\frac{1}{2}\partial_t\rho|\bm{u}|^2+\rho\partial_t\bm{u}\cdot\bm{u}\bigg)\d\x \\
        & =\int_\Omega\big(\nabla\cdot\big(2\eta\D(\bm{u})-p\I\big)+\mu\nabla\phi+\bm{J}\times\bm{B}+\rho\bm{g}\big)\cdot\bm{u}\d\x \\
        & =-\int_\Omega{2}\eta|\D(\bm{u})|^2\d\x+\int_\Omega(\mu\nabla\phi+\bm{J}\times\bm{B}+\rho\bm{g})\cdot\bm{u}\d\x.
    \end{aligned}
\end{equation*}
Next, using the definition of the chemical potential and the Cahn--Hilliard equation, the time derivative of the Helmholtz free energy satisfies
\begin{equation*}
    \begin{aligned}
        \frac{\d}{\d{t}}\int_\Omega\Bigg(\frac{\widehat{\lambda}\epsilon}{2}|\nabla\phi|^2+\frac{\widehat{\lambda}\epsilon^{-1}}{4}(\phi^2-1)^2\Bigg)\d\x & = \int_\Omega\Big(\widehat{\lambda}\epsilon\nabla\phi\cdot\nabla\partial_t\phi+\widehat{\lambda}\epsilon^{-1}(\phi^3-\phi)\partial_t\phi\Big)\d\x \\
        & =\int_\Omega\Big(-\widehat{\lambda}\epsilon\Delta\phi+\widehat{\lambda}\epsilon^{-1}(\phi^3-\phi)\Big)\partial_t\phi\d\x \\
        & =\int_\Omega\mu\partial_t\phi\d\x \\
        & =-\int_\Omega{m}_\epsilon|\nabla\mu|^2\d\x-\int_\Omega\mu\bm{u}\cdot\nabla\phi\d\x.
    \end{aligned}
\end{equation*}
Finally, using Faraday's law, Amp\`{e}re's law, and Ohm's law, we obtain
\begin{equation*}
    \begin{aligned}
        \frac{\d}{\d{t}}\int_\Omega\frac{1}{2}\mu_0^{-1}|\bm{B}|^2\d\x & =\int_\Omega\mu_0^{-1}\bm{B}\cdot\partial_t\bm{B}\d\x \\
        & =-\int_\Omega\mu_0^{-1}\bm{B}\cdot\nabla\times\bm{E}\d\x \\
        & =-\int_\Omega\mu_0^{-1}\nabla\times\bm{B}\cdot\bm{E}\d\x \\
        & =-\int_\Omega\sigma^{-1}|\bm{J}|^2\d\x+\int_\Omega\bm{J}\cdot\bm{u}\times\bm{B}\d\x.
    \end{aligned}
\end{equation*}
Summing the above three identities, we obtain the desired result.


\small
\bibliographystyle{abbrv}
\bibliography{ref}
\end{document}